\documentclass[preprint,12pt,nopreprintline]{elsarticle}
\usepackage[margin=1in]{geometry}  

\usepackage[T1]{fontenc}
\usepackage{lmodern}
\usepackage{amsmath,amssymb,amsthm}
\usepackage{graphicx}
\usepackage{enumitem}
\usepackage{booktabs}
\usepackage{xcolor}
\usepackage{microtype}
\usepackage[colorlinks=true,linkcolor=blue!50!black,citecolor=blue!50!black,urlcolor=blue!50!black]{hyperref}
\usepackage{thmtools}
\usepackage[nameinlink,capitalise,noabbrev]{cleveref}
\graphicspath{{figures/}}

\newtheorem{theorem}{Theorem}[section]
\newtheorem{proposition}[theorem]{Proposition}
\newtheorem{lemma}[theorem]{Lemma}
\newtheorem{corollary}[theorem]{Corollary}
\theoremstyle{definition}
\newtheorem{definition}[theorem]{Definition}
\newtheorem{example}[theorem]{Example}
\newtheorem{remark}[theorem]{Remark}
\numberwithin{equation}{section}

\crefname{corollary}{corollary}{corollaries}
\Crefname{corollary}{Corollary}{Corollaries}

\newcommand{\ARV}{\mathsf{ARV}}
\newcommand{\Corr}{\operatorname{Corr}}
\newcommand{\Cov}{\operatorname{Cov}}
\newcommand{\Var}{\operatorname{Var}}
\newcommand{\rhoo}{\rho_{11}^{\Lambda}}
\newcommand{\rhot}{\rho_{22}^{\Lambda}}
\newcommand{\E}{\mathbb{E}}
\newcommand{\R}{\mathbb{R}}
\newcommand{\1}{\mathbf 1}
\newcommand{\A}{\mathcal A}
\newcommand{\fold}[1]{#1^{\mathsf F}}  
\newcommand{\unfold}[1]{#1^{\mathsf U}}  

\begin{document}

\begin{frontmatter}

\title{The exact region of the Average Rectangular Volume and Spearman's rho, with an application to arachnitude}

\author[unisalento-econ]{Luca Anzilli\corref{cor1}}
\ead{luca.anzilli@unisalento.it}

\author[unisalento-math]{Fabrizio Durante}
\ead{fabrizio.durante@unisalento.it}

\cortext[cor1]{Corresponding author.}

\affiliation[unisalento-econ]{organization={Department of Economic Sciences, University of Salento},
            city={Lecce},
            country={Italy}}

\affiliation[unisalento-math]{organization={Department of Mathematics and Physics ``Ennio De Giorgi'', University of Salento},
            city={Lecce},
            country={Italy}}

\begin{abstract}
The average rectangular volume ($\ARV$) of a bivariate aggregation function
is the average of its volumes over all axis-parallel rectangles of the unit
square; for a copula, it is the expected copula volume of the rectangle
spanned by two independent uniform points. We determine the exact region of
the pairs $(\ARV(C),\rho(C))$, where $C$ ranges over all bivariate copulas
and $\rho$ is Spearman's rho. The region is compact, convex and symmetric,
and its upper boundary consists of three analytically parametrized arcs. A
previously proposed family of shuffles of Min generates one arc, but only on
part of its range; the other two arcs are generated by singular copulas
whose supports contain hyperbolic branches. The boundary is continuously
differentiable away from its left endpoint $\ARV=\frac{11}{120}$, where it
has a vertical tangent. The proof uses a folding transformation and
explicit Kantorovich potentials for a family of optimal transport problems.
Finally, we show that the arachnitude, a rank correlation used in
actuarial practice, is an affine transformation of $\ARV$, so that the
result also describes the joint region of Spearman's rho and arachnitude.
\end{abstract}


\end{frontmatter}

\section{Introduction}\label{sec:intro}

In the theory of aggregation functions, the volume that a bivariate
function assigns to a rectangle of the unit square measures the interaction
between its two arguments on that rectangle; it is non-negative on every
rectangle exactly when the function is supermodular (2-increasing), a
property shared by copulas and by the functions that generate universal
integrals \cite{GrabischMarichalMesiarPap2009,KlementMesiarPap2010,DuranteSempi2016}.
Motivated by this role, the \emph{average rectangular volume} ($\ARV$) of
an integrable function $F:[0,1]^2\to\R$ was introduced in
\cite{AnzilliDurante2026a} as the average of the $F$-volumes of all
axis-parallel rectangles of the unit square. For supermodular aggregation
functions it is a global index of the strength of the interaction between
the two inputs; for a copula $C$, it is the expected $C$-volume of the
rectangle spanned by two independent uniform points.

On the class of copulas, $\ARV$ behaves quite differently from the classical
measures of concordance. It is invariant under reflections of one argument,
so that it does not distinguish positive from negative dependence, and it
satisfies $\frac{11}{120}\le\ARV(C)\le\frac2{15}$
\cite[Proposition~3.8]{AnzilliDurante2026a}. It is therefore natural to ask
how much information $\ARV$ carries about a measure of concordance such as
Spearman's rho. 

For two functionals defined on copulas, this question is encoded in their
\emph{exact region}; determining such regions is often a challenging problem.
A classical example is the relation between Kendall's tau and Spearman's rho,
for which universal inequalities and sharp bounds have been established
\cite{Daniels1950,SchreyerPaulinTrutschnig2017}.
Related bounds and exact-region problems for restricted classes of copulas
and for other measures of association have been studied, for instance, in
\cite{KokolBukovsekEtAl2021,KokolBukovsekMojskerc2022,
KokolBukovsekStopar2023,KokolBukovsekStopar2024,
Tschimpke2025,TschimpkeSchreyerTrutschnig2025,
FuchsLimbachSchurrer2026,AnsariRockel2026,
KokolBukovsekEtAl2026,Rockel2026}.

Partial information on the $(\ARV,\rho)$ region was obtained in \cite{AnzilliDurante2026a,AnzilliDurante2026b}: the inequality
$|\rho(C)|\le3\,\ARV(C)+\frac35$ holds for every copula, minimal $\ARV$
forces $|\rho(C)|\le\frac12$, and a family of shuffles of Min traces two
curves joining the extreme parts of the region, which were conjectured to
be its boundary.

In this paper we give a complete description of the region
(\cref{thm:main}). Its upper boundary is the graph of a concave, strictly
increasing function $r^*$ and consists of three analytically parametrized
arcs. The shuffles of Min considered in  \cite{AnzilliDurante2026b} generate the
rightmost arc only for $0\le a\le\frac25$; for $\frac25<a<\frac12$ the
corresponding curve lies strictly inside the region. The remaining two arcs
are generated by singular copulas whose supports contain hyperbola-type
curves. The boundary is continuously differentiable away from its left
endpoint, and
\[
r^*(\alpha)=\frac12+4\left(\alpha-\frac{11}{120}\right)^{3/4}
+o\!\left(\left(\alpha-\frac{11}{120}\right)^{3/4}\right),
\qquad \alpha\downarrow\frac{11}{120},
\]
so that the boundary leaves $(\frac{11}{120},\frac12)$ with a vertical
tangent.

The $(\ARV,\rho)$ region is also of interest in practical applications. In
\cref{sec:arachnitude} we show that, for every copula $C$ and $(U,V)\sim C$,
\[
45\,\ARV(C)-5=\Corr\bigl((2U-1)^2,(2V-1)^2\bigr),
\]
i.e.\ $\ARV$ is an affine transformation of the \emph{arachnitude}. This
index is used in actuarial practice to calibrate the dependence between risk
drivers in internal capital models, since it measures whether two variables
tend to be simultaneously extreme, in either direction, a feature that
Spearman's rho cannot detect \cite{ShawSmithSpivak2011,AndroschuckEtAl2017};
it is the second diagonal coefficient of the Legendre expansion of the copula
density, Spearman's rho being the first \cite{McNeilNeslehovaSmith2025}.
Consequently, \cref{thm:main} describes exactly which combinations of
Spearman's rho and arachnitude are compatible with some copula, and thus
provides a consistency check when both coefficients are estimated from data
or specified by expert judgment. For instance, a Spearman's rho of $0.9$ is
compatible only with an arachnitude of at least $0.249$.

The paper is organized as follows. \Cref{sec:folding} recalls aggregation functions, copulas and $\ARV$,
introduces the folding transformation and establishes the basic geometry of
the region. \Cref{sec:families} constructs the three copula families.
\Cref{sec:main} states the main result and its consequences.
\Cref{sec:arachnitude} establishes the link with arachnitude and discusses
examples, while \Cref{sec:conclusions} concludes the paper.
The proof of the main result is deferred to \ref{app:proofs}.

Generative AI tools were used in the research process, for the numerical
exploration, for the construction of the extremal families and of the
Kantorovich potentials, and for symbolic verifications; all arguments have
been checked by the authors. Details are given in the declaration at the end
of the paper.

\section{Preliminary definitions and properties}
\label{sec:folding}

We first recall the basic notions on aggregation functions and copulas and
the definition of the average rectangular volume. A (bivariate) \emph{aggregation function} is a function
$A:[0,1]^2\to[0,1]$ that is non-decreasing in each argument and satisfies
$A(0,0)=0$ and $A(1,1)=1$ \cite{GrabischMarichalMesiarPap2009}. For a
function $F:[0,1]^2\to\R$ and a rectangle
\[
R=[x_1,y_1]\times[x_2,y_2],
\qquad
x_i\le y_i,
\]
the \emph{$F$-volume} of $R$ is
\[
V_F(R)
=
F(y_1,y_2)-F(x_1,y_2)-F(y_1,x_2)+F(x_1,x_2).
\]
$F$ is \emph{supermodular} (or \emph{$2$-increasing}) if $V_F(R)\ge0$ for
every rectangle $R$. A (bivariate) \emph{copula} is a supermodular
aggregation function with uniform margins, $C(u,1)=u$ and $C(1,v)=v$ for
$u,v\in[0,1]$ \cite{Nelsen2006,DuranteSempi2016}. Let $\mathcal C$ denote
the class of bivariate copulas and let
\[
M(u,v)=\min\{u,v\},
\qquad
W(u,v)=\max\{u+v-1,0\},
\qquad
\Pi(u,v)=uv .
\]
We identify a copula with the probability measure it induces on
$[0,1]^2$. We write $(U,V)\sim C$ to mean that $C$ is the joint
distribution function of the random pair $(U,V)$, i.e.\
$C(u,v)=\mathbb P(U\le u,V\le v)$; in particular, $U$ and $V$ are uniformly
distributed on $[0,1]$. In this case $V_C(R)=\mathbb P((U,V)\in R)$ for every
closed rectangle $R$, since the margins are continuous. Spearman's rho is
\[
\rho(C)
=
12\E[UV]-3
=
12\int_{[0,1]^2}C(u,v)\,du\,dv-3.
\]

For $\mathbf x\le\mathbf y$ in $[0,1]^2$ (componentwise), write
$[\mathbf x,\mathbf y]=[x_1,y_1]\times[x_2,y_2]$.

\begin{definition}[\cite{AnzilliDurante2026a}]\label{def:arv}
The \emph{average rectangular volume} of an integrable function
$F:[0,1]^2\to\R$ is
\[
\ARV(F)
=
\frac{1}{|\mathcal R_4|}
\int_{\mathcal R_4}
V_F([\mathbf x,\mathbf y])\,d\mathbf x\,d\mathbf y,
\]
where $\mathcal R_4=\{(\mathbf x,\mathbf y)\in[0,1]^2\times[0,1]^2:
\mathbf x\le\mathbf y\}$.
\end{definition}

By Fubini's theorem \cite{AnzilliDurante2026a},
\begin{equation}\label{eq:arv-lebesgue}
\ARV(F)=4\int_{[0,1]^2}F(u,v)(1-2u)(1-2v)\,du\,dv .
\end{equation}
For a copula, the normalized integral over $\mathcal R_4$ equals, by
symmetry, the expected $C$-volume of the rectangle $R_{\mathbf X,\mathbf Y}$
spanned by two independent uniform points $\mathbf X,\mathbf Y$ of
$[0,1]^2$ \cite[Remark~2.2]{AnzilliDurante2026a}. Since the interval
spanned by two independent uniform points of $[0,1]$ covers $u$ with
probability $2u(1-u)$, it holds that
\begin{equation}\label{eq:arv-expectation}
\ARV(C)=\mathbb P\bigl((U,V)\in R_{\mathbf X,\mathbf Y}\bigr)
=4\,\E[U(1-U)V(1-V)],
\end{equation}
where $(U,V)\sim C$ is independent of $(\mathbf X,\mathbf Y)$. On
$\mathcal C$, $\ARV$ is affine, invariant under the reflection of one
argument, and satisfies $\frac{11}{120}\le\ARV(C)\le\frac2{15}$
(see \cite{AnzilliDurante2026a}); under independence, $\ARV(\Pi)=\frac19$.

The main tool of this paper is a transformation that exploits the symmetry
of $\ARV$ with respect to the center of the unit square. Rectangles
centered at $(\frac12,\frac12)$ play a special role:
for $(u,v)\in[0,\frac12]^2$, let
\[
R_{u,v}
=
[u,1-u]\times[v,1-v]
\]
be the rectangle centered at $(1/2,1/2)$.

Let $(U,V)\sim C$ and define
\[
S=|2U-1|,
\qquad
T=|2V-1|.
\]
Since $U$ and $V$ are uniform on $[0,1]$, both $S$ and $T$ are again
uniform on $[0,1]$.

\begin{definition}
For a given copula $C$, the \emph{folded copula} of $C$, denoted by $\fold C$, is the copula of the random pair
$(S,T)=(|2U-1|,|2V-1|)$.
\end{definition}

Notice that this construction is related to the shuffle of copulas via suitable measure-preserving transformations of the unit interval (see \cite{DuranteSarkociSempi2009}).

By \emph{flipping one coordinate} we mean replacing $(U,V)$ by
$(1-U,V)$ or by $(U,1-V)$. Such a flip does not change the folded
copula; it preserves $\ARV$ and changes the sign of $\rho$.

For a copula $D$, we write $\E_D$ for expectation with respect to
$(S,T)\sim D$.

\begin{proposition}\label{prop:folding}
For every copula $C$ the following statements hold.
\begin{enumerate}[label=\textup{(\roman*)}]
\item
\[
V_C(R_{u,v})
=
\fold C(1-2u,1-2v),
\qquad
(u,v)\in\left[0,\frac12\right]^2.
\]

\item
\begin{equation}
\ARV(C)
=
\int_{[0,1]^2} st\,\fold C(s,t)\,ds\,dt
=
\frac1{12}
+
\frac14\,\E_{\fold C}[S^2T^2].
\label{eq:fold-arv}
\end{equation}

\item
\begin{equation}
|\rho(C)|
\le
3\,\E_{\fold C}[ST].
\label{eq:fold-rho}
\end{equation}
\end{enumerate}

Conversely, for every copula $D$ there exists a copula $C$ such that
\[
\fold C=D
\qquad\text{and}\qquad
\rho(C)=3\,\E_D[ST].
\]
\end{proposition}

\begin{proof}
Since the marginals of $C$ are continuous, $V_C(R_{u,v})=\mathbb P((U,V)\in R_{u,v})$ (boundaries carry no mass).
The event $(U,V)\in R_{u,v}$ is equivalent to
\[
S\le1-2u,
\qquad
T\le1-2v,
\]
because $u\le x\le1-u$ if and only if $|2x-1|\le1-2u$. Since $S$ and $T$ are uniform, the joint distribution function of
$(S,T)$ coincides with $\fold C$ on $[0,1]^2$, and \textup{(i)} follows.

To prove \textup{(ii)}, recall that the average over $\mathcal R_4$ in the definition of $\ARV$ equals, by symmetry, the
expected $C$-measure of the random rectangle $R_{\mathbf X,\mathbf Y}$ spanned by two independent uniform points
$\mathbf X,\mathbf Y$ of $[0,1]^2$ \cite[Remark~2.2]{AnzilliDurante2026a}. Fix a point $(x,y)$ in the unit square.
The probability that the random rectangle covers this point is
$4x(1-x)y(1-y)$ (each coordinate interval covers the corresponding coordinate with probability $2x(1-x)$, resp.\ $2y(1-y)$). It therefore holds that
\[
\ARV(C)=4\,\E[U(1-U)V(1-V)]
=\frac14\,\E[(1-S^2)(1-T^2)]
=\frac1{12}+\frac14\,\E[S^2T^2]
\]
(because $4U(1-U)=1-S^2$ and $4V(1-V)=1-T^2$).
Here we used $\E[S^2]=\E[T^2]=1/3$. Also, for any copula $D$,
\[
\int_{[0,1]^2}st\,D(s,t)\,ds\,dt
=\E_D\!\left[\int_S^1s\,ds\int_T^1t\,dt\right]
=\frac14\,\E_D[(1-S^2)(1-T^2)].
\]
This proves both equalities in \eqref{eq:fold-arv}.

For \textup{(iii)}, uniformity of the marginals gives
\[
\rho(C)=3\,\E[(2U-1)(2V-1)].
\]
Consequently,
\[
|\rho(C)|\le3\,\E[|2U-1|\,|2V-1|]
=3\,\E_{\fold C}[ST].
\]

For the converse, let $(S,T)\sim D$, let $\varepsilon$ be an independent
symmetric random sign, and set
\[
U=\frac{1+\varepsilon S}{2},
\qquad
V=\frac{1+\varepsilon T}{2}.
\]
Then $U$ and $V$ are uniform on $[0,1]$, the copula $C$ of $(U,V)$
satisfies $\fold C=D$, and
\[
\E[UV]
=
\frac14\bigl(1+\E_D[ST]\bigr).
\]
Therefore
\[
\rho(C)
=
12\E[UV]-3
=
3\,\E_D[ST].
\]
\end{proof}

We denote this particular \emph{unfolding} by $\unfold D$. Explicitly, conditioning on $\varepsilon$,
\[
\begin{split}\unfold D(u,v)=\tfrac12\Big[&D\big((2u-1)^+,(2v-1)^+\big)+D\big((1-2u)^+,(1-2v)^+\big)\\&+1-(1-2u)^+-(1-2v)^+\Big],\end{split}
\]
i.e.\ $\unfold D$ is the ordinal sum of the survival copula of $D$ on $[0,\frac12]^2$ and of $D$ on $[\frac12,1]^2$.
It uses a common independent sign in the two coordinates, so
\[
(\unfold D)^{\mathsf F}=D,\qquad
\rho(\unfold D)=3\E_D[ST].
\]
This construction will allow us to specify the extremal copulas on
the simpler folded square and then recover their original supports.

With these tools at hand, we can establish the basic geometric properties
of the region under study, namely the set of all pairs of values of the
average rectangular volume and of Spearman's rho that are attained by some
copula:
\[
\A=\{(\ARV(C),\rho(C)):C\in\mathcal C\}\subseteq\R^2 .
\]
Throughout, we denote by $(\alpha,r)$ the generic point of the plane, where
the first coordinate $\alpha$ refers to $\ARV$ and the second coordinate $r$
to Spearman's rho.
The class $\mathcal C$ is compact with respect to the uniform convergence. Both functionals are continuous and affine:
this is immediate from their representations as integrals of bounded
continuous functions. Thus $\A$ is compact and convex. Reflection in
one coordinate preserves $\ARV$ and reverses the sign of $\rho$,
so $\A$ is symmetric with respect to the horizontal axis
$\{(\alpha,r):r=0\}$, i.e.\ $(\alpha,r)\in\A$ if and only if $(\alpha,-r)\in\A$.

For any copula $D$ (in the role of a folded copula), the Fr\'echet-Hoeffding bound inequalities $W\le D\le M$ and
\eqref{eq:fold-arv} give, since $\frac1{12}+\frac14\E_D[S^2T^2]=\int st\,D(s,t)\,ds\,dt$ is monotone in $D$,
\[
\frac{11}{120}\le\frac1{12}+\frac14\E_D[S^2T^2]\le\frac2{15}.
\]
Indeed, $\E_W[S^2T^2]=\int_0^1s^2(1-s)^2\,ds=1/30$ and
$\E_M[S^2T^2]=\int_0^1s^4\,ds=1/5$. Both bounds are attained by
unfolding $W$ and $M$, respectively; convexity gives every intermediate
value. Moreover, equality at either endpoint forces $D=W$ or $D=M$:
the non-negative continuous function $D-W$ or $M-D$ has zero integral
against the strictly positive weight $st$ in the open square, and
hence vanishes everywhere by continuity.

Compactness and Proposition~\ref{prop:folding} now give
\begin{equation}
\begin{split}
r^*(\alpha)&:=\max\{\rho(C):\ARV(C)=\alpha\}\\
&=\max\left\{3\E_D[ST]:D\in\mathcal C,
\ \frac1{12}+\frac14\E_D[S^2T^2]=\alpha\right\},
\end{split}
\label{eq:rstar-folded}
\end{equation}
for $\alpha\in[11/120,2/15]$. Here ``$\le$'' in the second equality follows from Proposition~\ref{prop:folding}(ii)--(iii) applied to
$D=\fold C$, ``$\ge$'' from the unfolding $\unfold D$, and both maxima are attained because the constraint sets are closed
subsets of the compact set $\mathcal C$ and the functionals are continuous.

\begin{proposition}\label{prop:region-geometry}
The region is
\[
\A=\left\{(\alpha,r):\frac{11}{120}\le \alpha\le\frac2{15},\quad
|r|\le r^*(\alpha)\right\}.
\]
The function $r^*$ is continuous, concave and strictly increasing, with
\[
r^*(11/120)=\frac12,\qquad r^*(2/15)=1.
\]
\end{proposition}

\begin{proof}
Every vertical section of $\A$ over $[11/120,2/15]$ (the projection of $\A$ on the first axis, by the preceding
paragraph) is a nonempty compact interval,
symmetric about zero, which proves the displayed description.
Convexity gives concavity of $r^*$. Compactness gives upper
semicontinuity, and concavity gives continuity in the interior.
At each endpoint, comparison with a chord to a fixed interior point
gives the complementary lower limit inequality, so $r^*$ is continuous
on the closed interval as well.

The endpoint characterization above and \eqref{eq:rstar-folded} yield
\[
r^*(11/120)=3\int_0^1u(1-u)\,du=\frac12,\qquad
r^*(2/15)=3\int_0^1u^2\,du=1.
\]
Finally, $C\le M$ and the integral formula for $\rho$ imply
$\rho(C)\le1$, with equality only for $C=M$. Thus $r^*(\alpha)<1$ when
$\alpha<2/15$. For $\alpha_1<\alpha_2<2/15$, write
$\alpha_2=(1-\lambda)\alpha_1+\lambda(2/15)$ with $0<\lambda<1$.
Concavity gives
\[
r^*(\alpha_2)\ge(1-\lambda)r^*(\alpha_1)+\lambda>r^*(\alpha_1).
\]
The comparison with the right endpoint is immediate.
\end{proof}

\section{Three families of copulas}
\label{sec:families}

We now construct the copulas that will generate the upper boundary of the region $\mathcal{A}$. 
Family~I is a shuffle-of-Min family already considered in
\cite{AnzilliDurante2026b}; see also
\cite{MikusinskiSherwoodTaylor1992,DuranteSarkociSempi2009} for
background on shuffles. Families~II and~III are constructed on the
folded square and then unfolded. This makes their marginal constraints
and moment formulas transparent. Optimality will be proved in
\ref{app:proofs}.

\subsection{Family I: shuffles of Min}
\label{subsec:familyI}

For $a\in[0,\frac12]$, define
\[
C^{\mathrm I}_a(u,v)
=
\begin{cases}
	\max\{u+v-a,0\},
	& (u,v)\in[0,a]^2,\\[1mm]
	\max\{u+v-1,1-a\},
	& (u,v)\in[1-a,1]^2,\\[1mm]
	M(u,v),
	& \text{otherwise}.
\end{cases}
\]
Its flipped version is denoted by $\widetilde C^{\mathrm I}_a$ and is
given by
\[
\widetilde C^{\mathrm I}_a(u,v)
=
u-C^{\mathrm I}_a(u,1-v)
=
v-C^{\mathrm I}_a(1-u,v).
\]

The support of $C^{\mathrm I}_a$ is
\[
\bigl\{u+v=a\bigr\}\cap[0,a]^2
\;\cup\;
\bigl\{u=v\bigr\}\cap[a,1-a]^2
\;\cup\;
\bigl\{u+v=2-a\bigr\}\cap[1-a,1]^2.
\]
The central segment is omitted at the endpoint $a=1/2$. The support of $\widetilde C^{\mathrm I}_a$ is obtained by flipping one
coordinate.

\begin{figure}[htbp]
	\centering
	\includegraphics[width=.67\textwidth]{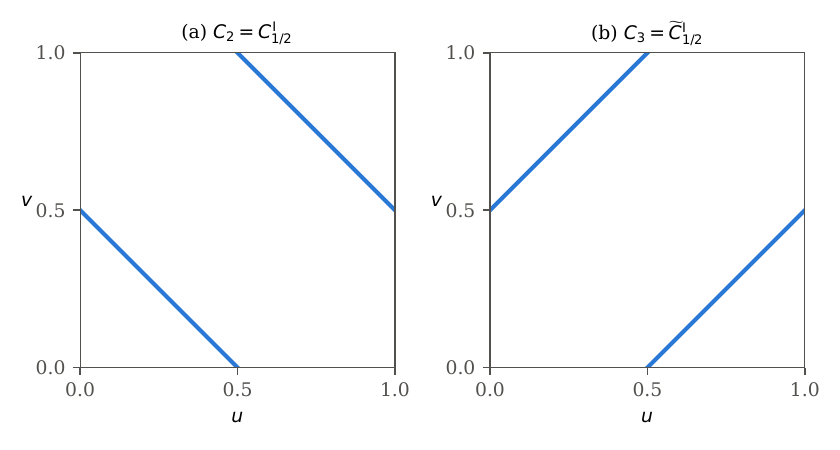}
	\caption{Supports of the endpoint copulas of Family~I:
		(a) $C_2=C^{\mathrm I}_{1/2}$ and
		(b) $C_3=\widetilde C^{\mathrm I}_{1/2}$.
		The support of $C_3$ is the reflection of that of $C_2$
		about $v=1/2$. Each segment carries probability mass $1/2$.}
	\label{fig:C2-C3-supports}
\end{figure}

In particular,
\[
C^{\mathrm I}_0=M,
\qquad
\widetilde C^{\mathrm I}_0=W,
\]
whereas, for later reference, we set
\[
C_2:=C^{\mathrm I}_{1/2},
\qquad
C_3:=\widetilde C^{\mathrm I}_{1/2}.
\]
The supports of $C_2$ and $C_3$ are shown in
Figure~\ref{fig:C2-C3-supports}.

A direct computation, detailed below, gives
\begin{equation}
	\ARV(C^{\mathrm I}_a)
	=
	\frac{2}{15}-\frac43(1-a)^2a^3,
	\qquad
	\rho(C^{\mathrm I}_a)
	=
	1-4a^3.
	\label{eq:familyI-arv-rho}
\end{equation}
Since $\ARV$ is invariant under flipping, while $\rho$ changes sign,
\[
\ARV(\widetilde C^{\mathrm I}_a)
=
\ARV(C^{\mathrm I}_a),
\qquad
\rho(\widetilde C^{\mathrm I}_a)
=
-\rho(C^{\mathrm I}_a).
\]

The folded representation of this family is particularly simple.
For $b\in[0,1]$, let $D_b$ be the copula of $(S,T_b(S))$, where $S$
is uniform and
\[
T_b(s)=\begin{cases}s,&0\le s\le b,\\1+b-s,&b<s\le1.\end{cases}
\]
Thus $D_b$ is the ordinal sum of a comonotone part on $[0,b]^2$
and a countermonotone part on $[b,1]^2$. Directly from the supports,
\begin{equation}
(C^{\mathrm I}_a)^{\mathsf F}=D_{1-2a},\qquad
C^{\mathrm I}_a=\unfold{D_{1-2a}}.
\label{eq:fold-familyI}
\end{equation}
In particular, $D_0=W$ and $D_1=M$.

The moments of Family~I are best computed on the folded square. By \eqref{eq:fold-familyI}, it suffices to compute
the moments of $D_b$. Since $D_b$ is concentrated on the diagonal of
$[0,b]^2$ and on the anti-diagonal of $[b,1]^2$,
\begin{align*}
\E_{D_b}[ST]
&=\frac{b^3}{3}+\int_b^1s(1+b-s)\,ds
=\frac{1+3b-3b^2+b^3}{6},\\
\E_{D_b}[S^2T^2]
&=\frac{b^5}{5}+\int_b^1s^2(1+b-s)^2\,ds
=\frac{1+5b+10b^2-10b^3-5b^4+5b^5}{30}.
\end{align*}
Substituting $b=1-2a$ in \eqref{eq:fold-arv} and in
$\rho(C^{\mathrm I}_a)=3\E_{D_b}[ST]$ gives
\eqref{eq:familyI-arv-rho}. Moreover,
\[
\frac{d}{db}\E_{D_b}[ST]=\frac{(1-b)^2}{2},
\qquad
\frac{d}{db}\E_{D_b}[S^2T^2]=\frac{(1-b)^2(1+b)(1+5b)}{6},
\]
which are the derivatives used in Remark~\ref{rem:kappa-familyI} and
in the proof of Proposition~\ref{prop:familyI-not-optimal}.

\subsection{Family II: a construction on the folded square}
\label{subsec:familyII}

Let $0<b\le1/5$ and $1/2+b\le c\le1$, and set
\[
\kappa=c(1+2b-c),\qquad b'=1+b-c.
\]
For $c\le t\le1$, define
\[
\xi_1(t)=\frac12\left(\frac\kappa t+t-1\right),\qquad
\xi_2(t)=\frac12\left(\frac\kappa t-t+1\right),
\]
and
\[
w_1(t)=\xi_1'(t)=\frac12\left(1-\frac\kappa{t^2}\right),\qquad
w_2(t)=-\xi_2'(t)=\frac12\left(1+\frac\kappa{t^2}\right).
\]
The parameter restrictions imply
\[
0\le b\le\kappa/2\le b'\le c,\qquad
w_1\ge0,\quad w_2>0,\quad w_1+w_2=1
\]
(see \eqref{eq:F2}--\eqref{eq:F3} below: $w_1\ge0$ is equivalent to $c^2\ge\kappa$, i.e.\ to $c\ge\frac12+b$).

The increasing branch $\xi_1$ joins $b$ to $\kappa/2$, whereas the
decreasing branch $\xi_2$ joins $b'$ to $\kappa/2$.

We define a measure $D_{b,c}$ through its integrals against bounded
Borel functions $f$:
\begin{align}
\int f\,dD_{b,c}
={}&\int_0^b f(s,s)\,ds
+\int_{b'}^c f(s,1+b-s)\,ds\nonumber\\
&+\sum_{i=1}^2\int_c^1w_i(t)
\bigl[f(\xi_i(t),t)+f(t,\xi_i(t))\bigr]\,dt.
\label{eq:familyII-measure}
\end{align}
This formula specifies both the support and the mass carried by each
branch; the support alone would not determine the measure.

\begin{figure}[t]
\centering
\includegraphics[width=\textwidth]{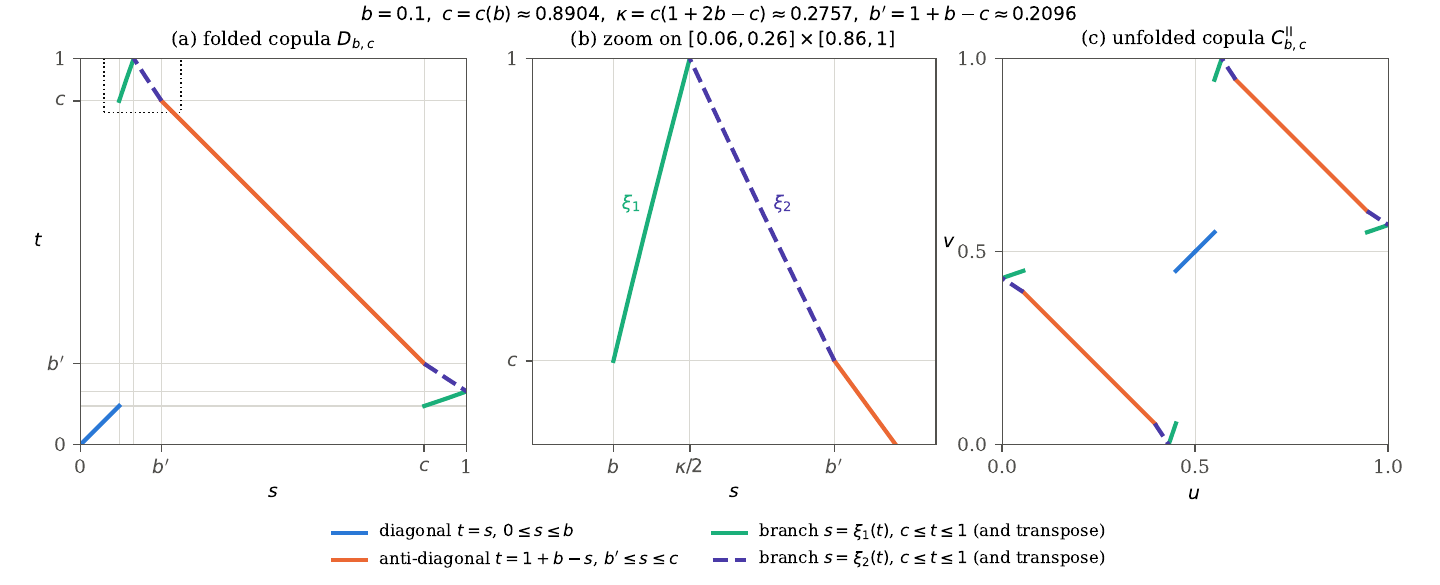}
\caption{The construction of Family~II for $b=0.1$ and $c=c(b)\approx0.8904$
(the parameters of the boundary arc $\Gamma^{\mathrm{II}}$; see
\cref{sec:main}). (a) Support of the folded copula $D_{b,c}$: comonotone
block, anti-diagonal and the two hyperbolic branches with their transposes.
(b) Zoom on the region where the pieces are glued together. (c) Support of
the unfolded copula $C^{\mathrm{II}}_{b,c}$.}
\label{fig:Dbc-support}
\end{figure}

\begin{remark}[On the construction of $D_{b,c}$]\label{rem:Dbc-construction}
The construction can be read as a rule that couples every value $t$ with
one or two values $s$ on the folded square; see \cref{fig:Dbc-support}.
\begin{enumerate}[label=\textup{(\alph*)}]
\item \emph{Three pieces.} The values $s\in[0,b]$ are coupled with
themselves (a comonotone block, as in Family~I). The values
$s\in[b',c]$ are coupled in reverse order along the anti-diagonal
$t=1+b-s$ (a countermonotone block). Finally, every $t$ in the top strip
$[c,1]$ is split between two partners $\xi_1(t)<\xi_2(t)$, with
probabilities $w_1(t)$ and $w_2(t)$; by symmetry, the same is done with the
roles of $s$ and $t$ exchanged.
\item \emph{Uniform marginals.} The comonotone and countermonotone blocks
use the $s$-values in $[0,b]\cup[b',c]$, and the transposed curved pieces
use $[c,1]$. The two branches must therefore cover the remaining interval
$[b,b']$ uniformly: $\xi_1$ maps $[c,1]$ onto $[b,\kappa/2]$, $\xi_2$ maps
$[c,1]$ onto $[\kappa/2,b']$, and their densities are $|\xi_i'|=w_i$. The
transposed copy has a uniform marginal on $[c,1]$ precisely because
$w_1+w_2=\xi_1'-\xi_2'=1$; since the two branches meet at $t=1$, this is
equivalent to $\xi_2(t)-\xi_1(t)=1-t$.
\item \emph{Hyperbolic branches.} The second relation,
$\xi_1(t)+\xi_2(t)=\kappa/t$, comes from optimality. For the transport
cost $c_\kappa(s,t)=st-s^2t^2/(2\kappa)$ introduced in \ref{app:proofs},
two distinct points $p,q$ can be coupled with the same $t$ by an optimal
plan only if
$p(1-pt/\kappa)=q(1-qt/\kappa)$, i.e.\ $(p+q)t=\kappa$ (see
\eqref{eq:slack-factorization}). The two relations give the formulas for
$\xi_1$ and $\xi_2$; each branch is a piece of a hyperbola, for instance
$s=\xi_1(t)$ satisfies $2st=t^2-t+\kappa$.
\item \emph{Gluing.} The choice $\kappa=c(1+2b-c)$ is exactly the condition
$\xi_1(c)=b$: the increasing branch starts where the comonotone block ends.
With this choice, $\xi_2(c)=1+b-c=b'$, so that the decreasing branch starts
at the upper end of the anti-diagonal. The condition $c\ge\frac12+b$ is
equivalent to $c^2\ge\kappa$, i.e.\ to $w_1\ge0$.
\end{enumerate}
\end{remark}

\begin{lemma}\label{lem:uniform-marginals}
The measure $D_{b,c}$ is a symmetric copula measure.
\end{lemma}
\begin{proof}
All terms in \eqref{eq:familyII-measure} are non-negative when $f\ge0$.
For the first marginal, the diagonal and anti-diagonal terms give
Lebesgue measure on $[0,b]$ and $[b',c]$, respectively.
Changing variables $s=\xi_1(t)$ and $s=\xi_2(t)$ in the untransposed
curved terms gives Lebesgue measure on $[b,\kappa/2]$ and
$[\kappa/2,b']$. The transposed terms give density $w_1+w_2=1$ on
$[c,1]$. These intervals partition $[0,1]$, up to endpoints, so the
first marginal is uniform and the total mass is one.
The measure is invariant under transposition: this is explicit for
the curved terms and follows from $s\mapsto1+b-s$ for the
anti-diagonal term. Its second marginal is therefore uniform as well.
When $c=1$, the curved integrals vanish and the same argument applies
with the corresponding empty intervals omitted.
\end{proof}

Define the copula of Family~II by
\[
C^{\mathrm{II}}_{b,c}:=\unfold{D_{b,c}}.
\]
Thus its measure is obtained by mapping $(s,t)$ to
$((1+s)/2,(1+t)/2)$ or $((1-s)/2,(1-t)/2)$, each with probability
$1/2$. The resulting copula is singular and radially symmetric.
For nondegenerate parameter values, its support consists of a central
diagonal segment, two anti-diagonal segments and eight curvilinear
arcs. The latter are images of hyperbolic branches: for example,
$s=\xi_1(t)$ satisfies $2st=t^2-t+\kappa$.
When $c=1$, the curved pieces disappear.

For clarity, write
\[
\mathcal D_b=\begin{cases}
\{(s,s):0\le s\le b\},&b>0,\\
\varnothing,&b=0.
\end{cases}
\]
The support on the folded square, including the extension to $b=0$
below, is
\begin{align}
\Sigma_{b,c}=\operatorname{supp}(D_{b,c})
={}&\mathcal D_b
\cup\{(s,1+b-s):b'\le s\le c\}\nonumber\\
&\cup\bigcup_{i=1}^2
\{(\xi_i(t),t),(t,\xi_i(t)):c\le t\le1\}.
\label{eq:folded-support}
\end{align}
In particular, when $b=0$ the origin is not an additional support point.

The moment calculation can be performed before unfolding. For $n=1,2$,
\begin{align}
\E_{D_{b,c}}[(ST)^n]
={}&\frac{b^{2n+1}}{2n+1}
+\int_{b'}^c[s(1+b-s)]^n\,ds\nonumber\\
&+2\sum_{i=1}^2\int_c^1 w_i(t)[t\xi_i(t)]^n\,dt.
\label{eq:familyII-moments}
\end{align}
The integrals in \eqref{eq:familyII-moments} are elementary. On the curved
branches, a direct computation gives
\[
\sum_{i=1}^2w_i(t)\,t\xi_i(t)=\frac{\kappa}{2t},
\qquad
\sum_{i=1}^2w_i(t)\,[t\xi_i(t)]^2
=\frac{t^2(1-t)^2}{4}+\frac{\kappa^2(2-t)}{4t},
\]
for $c\le t\le1$.
Hence \eqref{eq:familyII-moments} becomes
\begin{align*}
\E_{D_{b,c}}[ST]
&=\frac{b^3}{3}+\int_{b'}^cs(1+b-s)\,ds-\kappa\log c,\\
\E_{D_{b,c}}[S^2T^2]
&=\frac{b^5}{5}+\int_{b'}^cs^2(1+b-s)^2\,ds
+\frac12\int_c^1t^2(1-t)^2\,dt\\
&\qquad-\kappa^2\log c-\frac{\kappa^2(1-c)}{2}.
\end{align*}
With \eqref{eq:fold-arv}, $\rho(C^{\mathrm{II}}_{b,c})=3\E_{D_{b,c}}[ST]$
and $\kappa=c(1+2b-c)$, the remaining polynomial integrals give
\begin{align}
	\ARV\!\left(C^{\mathrm{II}}_{b,c}\right)
	={}&
	\frac{19}{240}
	+\frac{b^5}{24}
	-\frac{b^4}{24}
	-\frac{b^3}{12}
	-\frac{b^2}{12}
	-\frac{b}{24}
	\nonumber\\
	&\quad
	+\frac{2b^2c^3}{3}
	-\frac{b^2c^2}{2}
	-\frac{3bc^4}{4}
	+\frac{4bc^3}{3}
	-\frac{bc^2}{2}
	\nonumber\\
	&\quad
	+\frac{c^5}{5}
	-\frac{9c^4}{16}
	+\frac{c^3}{2}
	-\frac{c^2}{8}
	-\frac{c^2}{4}(c-1-2b)^2\log c,
	\label{eq:familyII-arv}\\[1mm]
	\rho\!\left(C^{\mathrm{II}}_{b,c}\right)
	={}&
	\frac{b^3}{2}
	-\frac{3b^2}{2}
	-\frac{3b}{2}
	+3bc^2
	-2c^3
	+3c^2
	-\frac12
	\nonumber\\
	&\quad
	+3c(c-1-2b)\log c.
	\label{eq:familyII-rho}
\end{align}

\subsection{Family III}
\label{subsec:familyIII}

The measure formula \eqref{eq:familyII-measure} remains valid at $b=0$:
the diagonal integral disappears, while the marginal calculation is
unchanged. Its integrals against continuous functions depend
continuously on the parameters, so this extension also agrees with
the weak limit of Family~II. For
$c\in[\frac12,1]$, we define
\begin{equation}
	C^{\mathrm{III}}_c
	:=
	C^{\mathrm{II}}_{0,c}.
	\label{eq:familyIII-definition}
\end{equation}
In this case
\[
\kappa=c(1-c),
\qquad
b'=1-c,
\]
and the central diagonal segment disappears. For $1/2<c<1$, its support consists of two anti-diagonal segments
and eight curvilinear arcs; at the endpoints, the degenerate pieces
are omitted.

For $c=1$, the curvilinear pieces disappear and
\[
C^{\mathrm{III}}_1=C_2=C^{\mathrm I}_{1/2}
\]
(here $\kappa=0$; the measure \eqref{eq:familyII-measure} is still well defined, since the curved integrals are over $[1,1]$,
and $D_{0,1}=W$).

Setting $b=0$ in \eqref{eq:familyII-arv}--\eqref{eq:familyII-rho} gives
\begin{align}
	\ARV\!\left(C^{\mathrm{III}}_c\right)
	&=
	\frac{19}{240}
	+\frac{c^5}{5}
	-\frac{9c^4}{16}
	+\frac{c^3}{2}
	-\frac{c^2}{8}
	-\frac{c^2(1-c)^2}{4}\log c,
	\label{eq:familyIII-arv}\\
	\rho\!\left(C^{\mathrm{III}}_c\right)
	&=
	-2c^3+3c^2-\frac12
	-3c(1-c)\log c.
	\label{eq:familyIII-rho}
\end{align}

Family~III is thus obtained from Family~II by removing the comonotone block:
on the folded square, small values of $S$ are no longer coupled with small
values of $T$, but are sent along the increasing branch $\xi_1$ to the top
strip $[c,1]$. In the original square this means that no mass is placed
close to the center along the main diagonal, which explains why these
copulas attain the smallest values of $\ARV$. As $c$ decreases from $1$ to
$\frac12$, the anti-diagonal piece $[1-c,c]$ shrinks and the curved branches
grow; for $c=1$ we recover $C_2$, whereas for $c=\frac12$ we have
$\kappa=\frac14$, the anti-diagonal piece reduces to a point, and the
support of $D_{0,1/2}$ consists only of the two branches and their
transposes. Only the range $c^*\le c\le1$ will be relevant for the
boundary.

At the other boundary of Family~II, namely $c=1$, we have
\begin{equation}
	C^{\mathrm{II}}_{b,1}
	=
	C^{\mathrm I}_{(1-b)/2}.
	\label{eq:familyII-familyI}
\end{equation}
At the folded level, these identities read
\begin{equation}
(C^{\mathrm{III}}_c)^{\mathsf F}=D_{0,c},\qquad D_{b,1}=D_b.
\label{eq:Dbc-Db}
\end{equation}
Thus Family~II provides a bridge between Families~III and~I.

\begin{figure}[t]
	\centering
	\includegraphics[width=\textwidth]{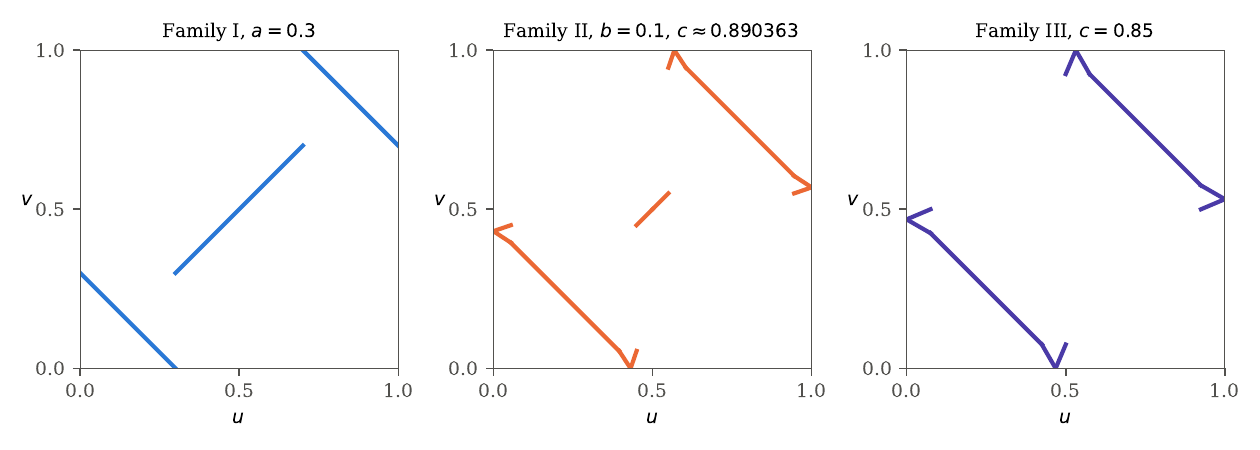}
	\caption{Representative supports of the three families.
		Left: Family~I with $a=0.30$.
		Center: Family~II with $b=0.10$ and $c\approx0.890363$.
		Right: Family~III with $c=0.85$.}
	\label{fig:three-families-support}
\end{figure}

Figure~\ref{fig:three-families-support} shows representative supports
of the three families.

\section{The exact boundary and the main result}
\label{sec:main}

We now identify the portions of the three families introduced in
Section~\ref{sec:families} that form the upper boundary of the
$(\ARV,\rho)$ region.

\subsection{The three boundary arcs}
\label{subsec:three-curves}

For Family~I, define
\begin{equation}
	\Gamma^{\mathrm I}
	=
	\left\{
	\left(
	\ARV(C^{\mathrm I}_a),
	\rho(C^{\mathrm I}_a)
	\right):
	0\le a\le\frac25
	\right\}.
	\label{eq:GammaI}
\end{equation}

For Family~II, introduce the polynomial
\begin{align}
	P(b,c)
	={}&
	6c^4-20c^3+18c^2-6c+1
	-8bc^3+36bc^2-24bc+4b
	\nonumber\\
	&\quad
	-6b^2c^2-18b^2c+6b^2
	+12b^3c+4b^3-5b^4.
	\label{eq:P}
\end{align}

Let $c^*$ be the unique root in $(1/2,1)$ of
\begin{equation}
	P(0,c)
	=
	6c^4-20c^3+18c^2-6c+1
	=
	0,
	\label{eq:cstar}
\end{equation}
so that
\[
c^*\approx0.7803501.
\]

For $b\in(0,\frac15]$, let $c(b)$ denote the unique solution of
\begin{equation}
	P(b,c)=0,
	\qquad
	\frac12+b<c\le1.
	\label{eq:cb}
\end{equation}
We set
\[
c(0):=c^*.
\]
As will be proved later, $c(b)$ is well defined on
$[0,\frac15]$, satisfies
\[
\lim_{b\downarrow0}c(b)=c^*,
\qquad
c\!\left(\frac15\right)=1,
\]
and selects the relevant one-dimensional subfamily of Family~II.

We then define
\begin{equation}
	\Gamma^{\mathrm{II}}
	=
	\left\{
	\left(
	\ARV(C^{\mathrm{II}}_{b,c(b)}),
	\rho(C^{\mathrm{II}}_{b,c(b)})
	\right):
	0\le b\le\frac15
	\right\},
	\label{eq:GammaII}
\end{equation}
where at $b=0$ we use the extension
$C^{\mathrm{II}}_{0,c}=C^{\mathrm{III}}_c$ introduced in
Section~\ref{subsec:familyIII}.

Finally, for Family~III, define
\begin{equation}
	\Gamma^{\mathrm{III}}
	=
	\left\{
	\left(
	\ARV(C^{\mathrm{III}}_c),
	\rho(C^{\mathrm{III}}_c)
	\right):
	c^*\le c\le1
	\right\}.
	\label{eq:GammaIII}
\end{equation}

The three arcs join continuously. The arcs
$\Gamma^{\mathrm{III}}$ and $\Gamma^{\mathrm{II}}$ meet at
\[
\left(
\ARV(C^{\mathrm{III}}_{c^*}),
\rho(C^{\mathrm{III}}_{c^*})
\right)
\approx
(0.0917551,0.5039872).
\]
The arcs $\Gamma^{\mathrm{II}}$ and $\Gamma^{\mathrm I}$ meet at
\[
\left(
\frac{962}{9375},
\frac{93}{125}
\right),
\]
corresponding to
\[
b=\frac15,
\qquad
c=1,
\qquad
a=\frac25.
\]

Thus
\[
\Gamma^{\mathrm{III}},
\qquad
\Gamma^{\mathrm{II}},
\qquad
\Gamma^{\mathrm I}
\]
join to form a continuous curve connecting
\[
\left(
\frac{11}{120},
\frac12
\right)
\qquad\text{and}\qquad
\left(
\frac{2}{15},
1
\right).
\]

\begin{figure}[t]
	\centering
	\includegraphics[width=.8\textwidth]{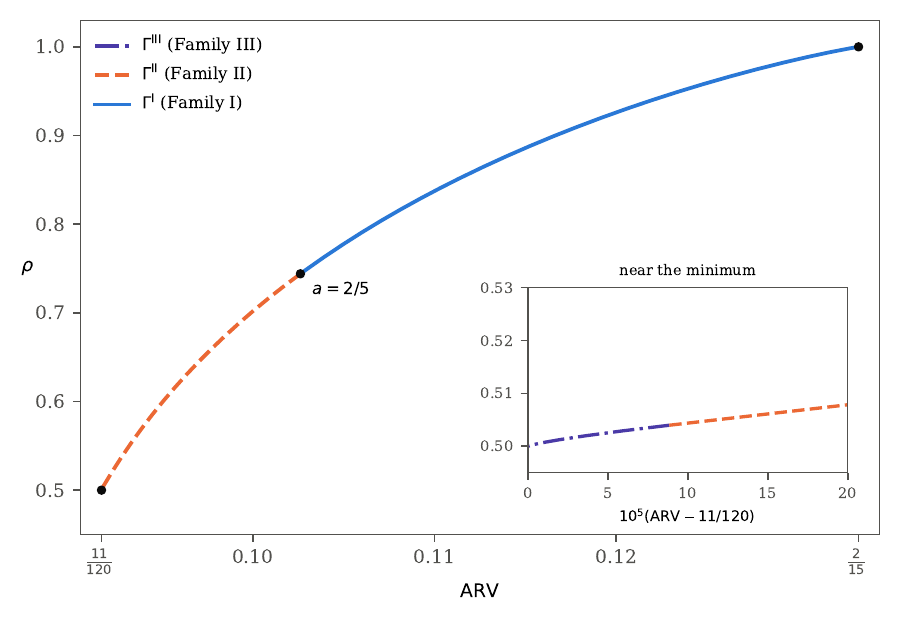}
	\caption{The three boundary arcs
		$\Gamma^{\mathrm{III}}$, $\Gamma^{\mathrm{II}}$ and
		$\Gamma^{\mathrm I}$ in the $(\ARV,\rho)$-plane.
        The inset resolves the left endpoint and the junction $c=c^*$;
        its horizontal coordinate is $10^5(\ARV-11/120)$.
        Distinct line styles identify the arcs in monochrome.}
	\label{fig:three-arcs}
\end{figure}

Figure~\ref{fig:three-arcs} shows the three arcs in the
$(\ARV,\rho)$-plane.

\subsection{Main result}
\label{subsec:main-result}

We can now state the main result.

\begin{theorem}\label{thm:main}
	Let
	\[
	\A
	=
	\left\{
	(\ARV(C),\rho(C)):
	C\in\mathcal C
	\right\}
	\]
	be the $(\ARV,\rho)$ region of bivariate copulas. Then
	\[
	\A
	=
	\left\{
	(\alpha,r)\in\mathbb R^2:
	\frac{11}{120}\le \alpha\le\frac{2}{15},
	\quad
	|r|\le r^*(\alpha)
	\right\},
	\]
	where $r^*$ is concave and strictly increasing and
	\[
	\operatorname{graph}(r^*)
	=
	\Gamma^{\mathrm{III}}
	\cup
	\Gamma^{\mathrm{II}}
	\cup
	\Gamma^{\mathrm I}.
	\]
	
    The arcs are generated by the following copulas (in order from left to right):
    \[
    \begin{array}{c|c|c}
    \text{arc}&\text{copula}&\text{parameter range}\\ \hline
    \Gamma^{\mathrm{III}}&C^{\mathrm{III}}_c&c^*\le c\le1\\
    \Gamma^{\mathrm{II}}&C^{\mathrm{II}}_{b,c(b)}&0\le b\le1/5\\
    \Gamma^{\mathrm I}&C^{\mathrm I}_a&0\le a\le2/5
    \end{array}
    \]
	The lower boundary is obtained by flipping one coordinate, which
	preserves $\ARV$ and changes the sign of $\rho$.
\end{theorem}

\begin{figure}[!p]
\centering
\includegraphics[width=\textwidth]{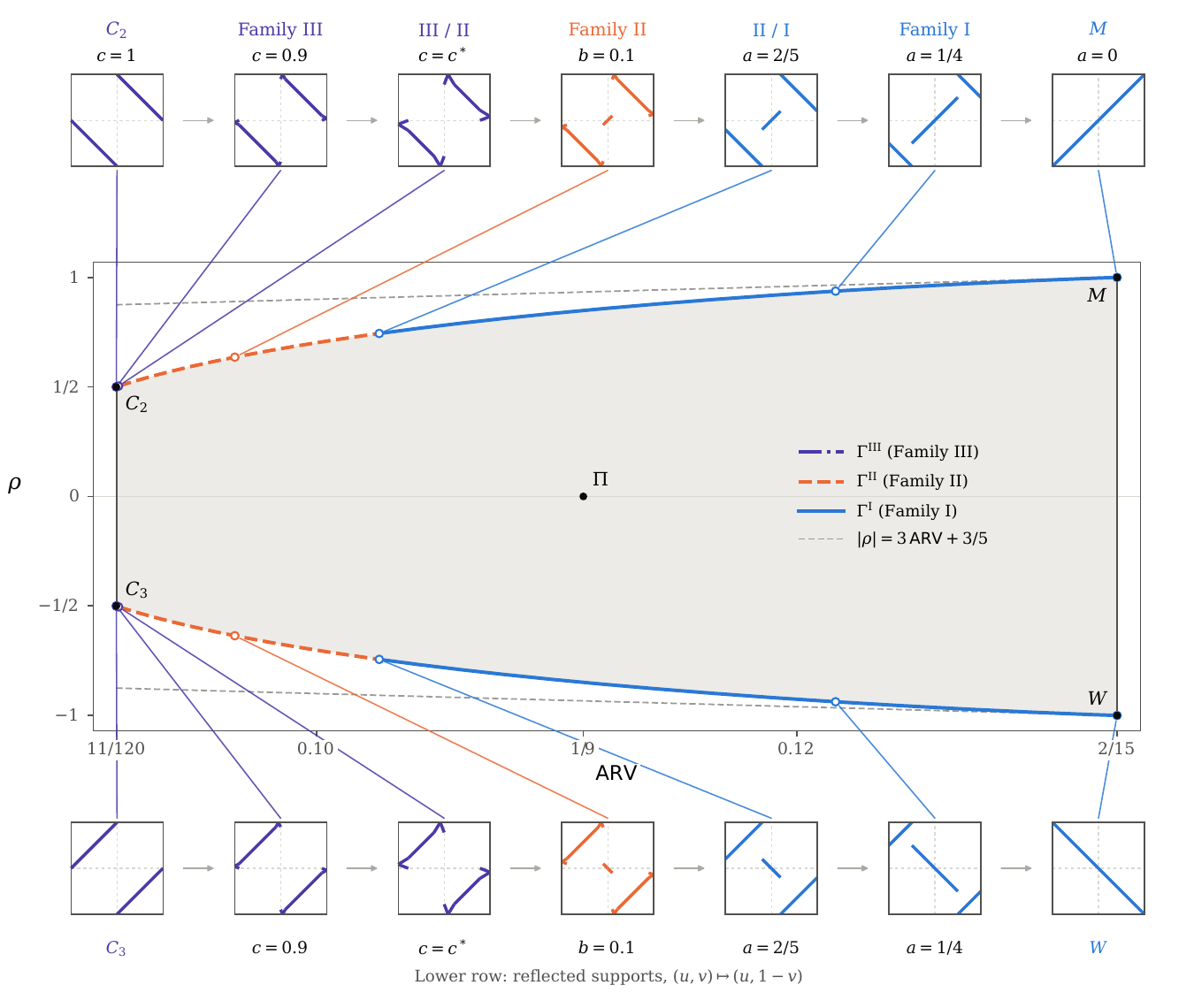}
\caption{The $(\ARV,\rho)$ region $\A$ (shaded). The upper and lower
boundaries are formed by the arcs $\Gamma^{\mathrm{III}}$,
$\Gamma^{\mathrm{II}}$ and $\Gamma^{\mathrm I}$ and by their flipped
versions; the vertical segment joins $C_2$ and $C_3$. Open circles mark the
junctions of the arcs. The top row shows the supports of the copulas
$C_2$, $C^{\mathrm{III}}_{0.9}$, $C^{\mathrm{III}}_{c^*}$,
$C^{\mathrm{II}}_{0.1,c(0.1)}$, $C^{\mathrm I}_{2/5}$,
$C^{\mathrm I}_{1/4}$ and $M$, which attain the upper boundary; the bottom
row shows the reflected supports, $(u,v)\mapsto(u,1-v)$, of the copulas
attaining the lower boundary, from $C_3$ to $W$. Since the points of
$\Gamma^{\mathrm{III}}$ are within $10^{-4}$ of $C_2$ in the $\ARV$
coordinate, they are not resolved at this scale (see
\cref{fig:region-zoom}). The dashed lines show the earlier bound
$|\rho|\le3\,\ARV+\frac35$ of \cite{AnzilliDurante2026a}.}
\label{fig:region-full}
\end{figure}

\begin{figure}[!htb]
\centering
\includegraphics[width=.75\textwidth]{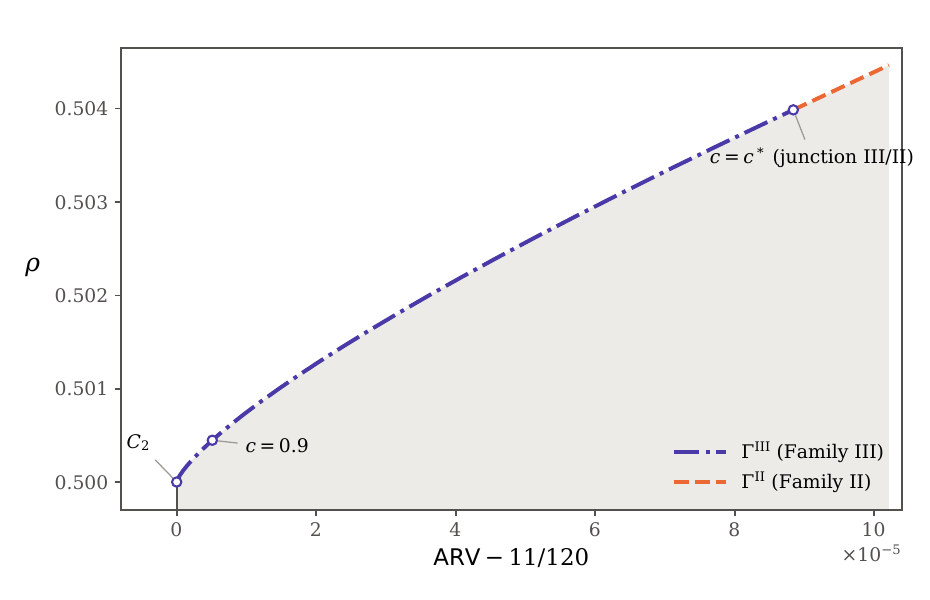}
\caption{Enlargement of the upper boundary of $\A$ near its left endpoint
$C_2=(\frac{11}{120},\frac12)$; the horizontal coordinate is
$\ARV-\frac{11}{120}$. The arc $\Gamma^{\mathrm{III}}$ leaves $C_2$ with a
vertical tangent (\cref{cor:vertical-tangent}) and meets
$\Gamma^{\mathrm{II}}$ at $c=c^*$.}
\label{fig:region-zoom}
\end{figure}

The region is represented in \cref{fig:region-full}, together with the
supports of some of the copulas attaining its boundary; an enlargement near
the left endpoint $C_2$ is shown in \cref{fig:region-zoom}. \Cref{thm:main}
can be stated as a sharp inequality between the two functionals.

\begin{corollary}\label{cor:sharp-arv}
For every copula $C$,
\[
|\rho(C)|\le r^*\bigl(\ARV(C)\bigr),
\]
and the inequality is sharp for every value of $\ARV(C)$. Equivalently,
$\ARV(C)\ge(r^*)^{-1}(|\rho(C)|)$ whenever $|\rho(C)|\ge\frac12$, while
every value $\ARV\in[\frac{11}{120},\frac2{15}]$ is compatible with
$|\rho|\le\frac12$.
\end{corollary}

The boundary is smooth away from its left endpoint, where it has a vertical
tangent. The proofs of the following two results are given in
\cref{subsec:vertical-tangent}.

\begin{corollary}\label{cor:C1}
The function $r^*$ is continuously differentiable on
$(11/120,2/15]$, with a one-sided derivative at $2/15$. Its derivative is
determined by the parameter $\kappa$ of the supporting cost:
\begin{align}
(r^*)'\bigl(\ARV(C^{\mathrm I}_a)\bigr)
&=\frac{6}{\kappa_{1-2a}}=\frac{9}{(1-a)(3-5a)},
&&0\le a\le\tfrac25,\nonumber\\
(r^*)'\bigl(\ARV(C^{\mathrm{II}}_{b,c(b)})\bigr)
&=\frac{6}{c(b)\bigl(1+2b-c(b)\bigr)},
&&0\le b\le\tfrac15,\label{eq:rstar-derivative}\\
(r^*)'\bigl(\ARV(C^{\mathrm{III}}_c)\bigr)
&=\frac{6}{c(1-c)},
&&c^*\le c<1.\nonumber
\end{align}
In particular, $(r^*)'(2/15)=3$, the three arcs join without corners,
with slopes $15$ at $\Gamma^{\mathrm{II}}\cap\Gamma^{\mathrm I}$ and
$6/(c^*(1-c^*))\approx35.005$ at
$\Gamma^{\mathrm{III}}\cap\Gamma^{\mathrm{II}}$, and
$(r^*)'(\alpha)\to+\infty$ as $\alpha\downarrow11/120$.
\end{corollary}

\begin{corollary}\label{cor:vertical-tangent}
As $\alpha\downarrow11/120$,
\[
r^*(\alpha)
=
\frac12
+
4\left(\alpha-\frac{11}{120}\right)^{3/4}
+
o\!\left(
\left(\alpha-\frac{11}{120}\right)^{3/4}
\right).
\]
In particular, the upper boundary has a vertical tangent at
$\bigl(11/120,1/2\bigr)$; by symmetry, the lower boundary has a
vertical tangent at $\bigl(11/120,-1/2\bigr)$.
\end{corollary}

The proof of \cref{thm:main}, together with those of
\Cref{cor:C1,cor:vertical-tangent}, is given in \ref{app:proofs}.
It relies on a reformulation of the boundary problem as an optimal
transport problem on the folded square and on the construction of explicit
Kantorovich potentials for the three extremal families.

\section{Average rectangular volume and arachnitude}
\label{sec:arachnitude}

In this section we show that, on copulas, $\ARV$ is an affine function of
the arachnitude, a second-order rank correlation used in actuarial
practice, and we restate \cref{thm:main} as the exact joint region of
Spearman's rho and arachnitude. Spearman's rho measures the \emph{direction}
of dependence; arachnitude measures its \emph{magnitude}, i.e.\ whether the
two ranks tend to be simultaneously extreme, in either tail, or
simultaneously central. It was introduced to calibrate the dependence
between risk drivers in internal capital models
\cite{ShawSmithSpivak2011,AndroschuckEtAl2017}, and is a prominent example of
the non-monotonic rank correlations studied in
\cite{McNeilNeslehovaSmith2025}.

We first recall the Legendre rank correlations.

Let $\mathcal C$ be the class of bivariate copulas and $(U,V)\sim C$. For
non-constant $g,h\in L^2(0,1)$, the generalized Spearman correlation of
\cite{McNeilNeslehovaSmith2025} is $\Corr(g(U),h(V))$; if $g$ and $h$ are
centered and normalized in $L^2(0,1)$, it equals $\E[g(U)h(V)]$. Let $P_j$
be the Legendre polynomial of degree $j$ on $[-1,1]$ and
\[
\Lambda_j(u)=\sqrt{2j+1}\,P_j(2u-1),\qquad j\ge0,
\]
the orthonormal shifted Legendre basis of $L^2(0,1)$. The
\emph{Legendre rank correlations} of $C$ are
$\rho_{jk}^{\Lambda}(C)=\E[\Lambda_j(U)\Lambda_k(V)]$, $j,k\ge1$; if $C$ has a
square-integrable density $c$, then
$c(u,v)=1+\sum_{j,k\ge1}\rho_{jk}^{\Lambda}(C)\Lambda_j(u)\Lambda_k(v)$ in
$L^2$. Since
\begin{equation}\label{eq:L1L2}
\Lambda_1(u)=\sqrt3\,(2u-1),
\qquad
\Lambda_2(u)=\sqrt5\,(6u^2-6u+1)=\frac{\sqrt5}{2}\bigl(3(2u-1)^2-1\bigr),
\end{equation}
we have
\begin{equation}\label{eq:rho11}
\rho_{11}^{\Lambda}(C)=3\,\E[(2U-1)(2V-1)]=12\,\E[UV]-3=\rho(C),
\end{equation}
Spearman's rho, and, since $\Lambda_2$ is the standardized version of
$(2u-1)^2$,
\begin{equation}\label{eq:rho22}
\rho_{22}^{\Lambda}(C)=\Corr\bigl((2U-1)^2,(2V-1)^2\bigr),
\end{equation}
the arachnitude \cite{ShawSmithSpivak2011,McNeilNeslehovaSmith2025}. A flip
$u\mapsto1-u$ changes the sign of $\Lambda_1$ and leaves $\Lambda_2$
unchanged; hence flipping one coordinate changes the sign of
$\rho_{11}^{\Lambda}$ and preserves $\rho_{22}^{\Lambda}$. In particular
$\rho_{22}^{\Lambda}(M)=\rho_{22}^{\Lambda}(W)=1$.

The second of these coefficients is an affine function of the average
rectangular volume.

\begin{theorem}\label{thm:arv-arachnitude}
For every bivariate copula $C$,
\begin{equation}\label{eq:arv-rho22}
\rho_{22}^{\Lambda}(C)=45\,\ARV(C)-5,
\qquad\text{equivalently}\qquad
\ARV(C)=\frac19+\frac1{45}\,\rho_{22}^{\Lambda}(C).
\end{equation}
\end{theorem}

\begin{proof}
By \eqref{eq:L1L2}, $u(1-u)=\frac16-\frac1{6\sqrt5}\Lambda_2(u)$. Since
$\E[\Lambda_2(U)]=\E[\Lambda_2(V)]=0$, \eqref{eq:arv-expectation} gives
\[
\ARV(C)=4\,\E\Bigl[\Bigl(\frac16-\frac{\Lambda_2(U)}{6\sqrt5}\Bigr)
\Bigl(\frac16-\frac{\Lambda_2(V)}{6\sqrt5}\Bigr)\Bigr]
=\frac19+\frac1{45}\,\E[\Lambda_2(U)\Lambda_2(V)].
\qedhere
\]
\end{proof}

Under independence, $\ARV(\Pi)=\frac19$ and $\rho_{22}^{\Lambda}(\Pi)=0$.
Positive arachnitude means that random rectangles capture \emph{more} copula
mass than under independence, i.e.\ that $U(1-U)$ and $V(1-V)$ are
positively correlated: the two ranks tend to be simultaneously central or
simultaneously extreme. Negative arachnitude means that one rank tends to be
extreme when the other is central.

Spearman's rho is also an expected volume of random rectangles, but of
rectangles attached to the corners of the square. For $\mathbf x\in[0,1]^2$
and a corner $\mathbf e\in\{0,1\}^2$, let $R_{\mathbf e}(\mathbf x)$ be the
rectangle with opposite vertices $\mathbf x$ and $\mathbf e$.

\begin{proposition}\label{prop:rho-rectangles}
Let $\mathbf X$ be uniform on $[0,1]^2$. For every copula $C$,
$$
\rho_{11}^{\Lambda}(C)
=3\,\E\Bigl[V_C\bigl(R_{(0,0)}(\mathbf X)\bigr)+V_C\bigl(R_{(1,1)}(\mathbf X)\bigr)
-V_C\bigl(R_{(1,0)}(\mathbf X)\bigr)-V_C\bigl(R_{(0,1)}(\mathbf X)\bigr)\Bigr].
$$
\end{proposition}

\begin{proof}
The four rectangles cover $(u,v)$ with probabilities $(1-u)(1-v)$, $uv$,
$u(1-v)$ and $(1-u)v$, whose signed sum is $(1-2u)(1-2v)$. As in
\eqref{eq:arv-expectation}, the right-hand side equals
$3\,\E[(1-2U)(1-2V)]=\rho_{11}^{\Lambda}(C)$.
\end{proof}

A flip exchanges the corners $(0,0)\leftrightarrow(1,0)$ and
$(0,1)\leftrightarrow(1,1)$, so the anchored representation changes sign,
whereas the law of the spanned rectangle $R_{\mathbf X,\mathbf Y}$ is
invariant. Spearman's rho measures the direction of dependence because its
rectangles are anchored; arachnitude measures its magnitude because its
rectangles are attached to nothing.

As a consequence, every result on $\ARV$ can be translated into a result on
arachnitude.

By \cref{thm:arv-arachnitude}, the map
\begin{equation}\label{eq:translation}
\Phi(\alpha,r)=(45\alpha-5,\,r)
\end{equation}
is an increasing affine bijection of the plane that carries the
$(\ARV,\rho)$ region onto the joint region of arachnitude and Spearman's rho:
\begin{equation}\label{eq:legendre-region}
\A^{\Lambda}:=\bigl\{(\rho_{22}^{\Lambda}(C),\rho_{11}^{\Lambda}(C)):C\in\mathcal C\bigr\}
=\Phi(\A).
\end{equation}
Boundaries are mapped onto boundaries, slopes are divided by $45$, and the
first coordinate is translated according to the following dictionary.

\begin{table}[ht]
\centering
\caption{Correspondence between $\ARV$ and arachnitude,
$\rho_{22}^{\Lambda}=45\,\ARV-5$.}
\label{tab:dictionary}
\begin{tabular}{lcc}
\toprule
& $\ARV$ & $\rho_{22}^{\Lambda}$\\
\midrule
Minimum (e.g.\ $C_2$, $C_3$) & $\frac{11}{120}$ & $-\frac78$\\[1mm]
Independence $\Pi$ & $\frac19$ & $0$\\[1mm]
Maximum ($M$, $W$, $\frac{M+W}2$) & $\frac2{15}$ & $1$\\[1mm]
Junction of $\Gamma^{\mathrm{III}}$ and $\Gamma^{\mathrm{II}}$ & $0.0917551$ & $-0.871022$\\[1mm]
Junction of $\Gamma^{\mathrm{II}}$ and $\Gamma^{\mathrm I}$ & $\frac{962}{9375}$ & $-\frac{239}{625}$\\[1mm]
Slope of a supporting line & $\frac6\kappa$ & $\frac2{15\kappa}$\\
\bottomrule
\end{tabular}
\end{table}

In particular, the range $\frac{11}{120}\le\ARV\le\frac2{15}$ of
\cite{AnzilliDurante2026a} gives the sharp range of arachnitude,
\[
-\frac78\le\rho_{22}^{\Lambda}(C)\le1,
\]
a value also obtained in \cite{McNeilNeslehovaSmith2025}. By the
discussion preceding \cref{prop:region-geometry}, the lower bound is attained
exactly by the copulas whose folded copula is $W$ and the upper bound exactly
by those whose folded copula is $M$.

In particular, this applies to the main result of this paper.

By \cref{thm:arv-arachnitude}, \cref{thm:main} and
\cref{cor:C1,cor:vertical-tangent} translate into statements on the joint region
$\A^{\Lambda}$ of \eqref{eq:legendre-region}. For $\eta\in[-\frac78,1]$ put
\begin{equation}\label{eq:Rstar}
R^*(\eta)=r^*\Bigl(\frac{\eta+5}{45}\Bigr)
=\max\bigl\{\rho_{11}^{\Lambda}(C):\ C\in\mathcal C,\ \rho_{22}^{\Lambda}(C)=\eta\bigr\}.
\end{equation}

\begin{theorem}\label{thm:main-legendre}
The joint region of arachnitude and Spearman's rho is
\[
\A^{\Lambda}
=
\left\{
(\eta,r)\in\R^2:\ -\frac78\le\eta\le1,\quad |r|\le R^*(\eta)
\right\},
\]
where $R^*$ is continuous, concave and strictly increasing, with
$R^*(-\frac78)=\frac12$ and $R^*(1)=1$. The graph of $R^*$ is the union of
the three arcs $\Phi(\Gamma^{\mathrm{III}})$, $\Phi(\Gamma^{\mathrm{II}})$
and $\Phi(\Gamma^{\mathrm I})$, generated by the copulas of
\cref{thm:main}:
\begin{enumerate}[label=\textup{(\alph*)}]
\item for $c^*\le c\le1$, the points
$\bigl(\rho_{22}^{\Lambda}(C^{\mathrm{III}}_c),\rho_{11}^{\Lambda}(C^{\mathrm{III}}_c)\bigr)$,
with
\[
\rho_{22}^{\Lambda}(C^{\mathrm{III}}_c)
=-\frac{23}{16}+9c^5-\frac{405c^4}{16}+\frac{45c^3}{2}-\frac{45c^2}{8}
-\frac{45c^2(1-c)^2}{4}\log c,
\]
and $\rho_{11}^{\Lambda}(C^{\mathrm{III}}_c)=\rho(C^{\mathrm{III}}_c)$ given
by \eqref{eq:familyIII-rho};
\item for $0\le b\le\frac15$, the points
$\bigl(\rho_{22}^{\Lambda}(C^{\mathrm{II}}_{b,c(b)}),\rho_{11}^{\Lambda}(C^{\mathrm{II}}_{b,c(b)})\bigr)$,
with $\rho_{22}^{\Lambda}=45\,\ARV-5$ and $\ARV$, $\rho$ given by
\eqref{eq:familyII-arv}--\eqref{eq:familyII-rho};
\item for $0\le a\le\frac25$, the points
\[
\bigl(\rho_{22}^{\Lambda}(C^{\mathrm I}_a),\rho_{11}^{\Lambda}(C^{\mathrm I}_a)\bigr)
=\bigl(1-60(1-a)^2a^3,\ 1-4a^3\bigr).
\]
\end{enumerate}
The arcs meet at $(\rho_{22}^{\Lambda},\rho_{11}^{\Lambda})\approx(-0.871022,0.503987)$
and at $(-\frac{239}{625},\frac{93}{125})$, and join $(-\frac78,\frac12)$ to
$(1,1)$. The lower boundary is obtained by flipping one coordinate, which
preserves $\rho_{22}^{\Lambda}$ and changes the sign of
$\rho_{11}^{\Lambda}$.
\end{theorem}

\begin{proof}
By \eqref{eq:legendre-region}, $\A^{\Lambda}=\Phi(\A)$, and since $\Phi$
only rescales and translates the first coordinate, \cref{thm:main} gives the
description of $\A^{\Lambda}$ with $R^*=r^*\circ\Phi_1^{-1}$, where
$\Phi_1(\alpha)=45\alpha-5$; continuity, concavity and strict monotonicity
are preserved. The endpoint values follow from
\cref{prop:region-geometry}. The parametrizations follow from
\eqref{eq:familyI-arv-rho}, \eqref{eq:familyIII-arv} and
\eqref{eq:familyII-arv} by applying $\Phi_1$, and the junction points from
those of \cref{subsec:three-curves}.
\end{proof}

\begin{corollary}\label{cor:sharp-legendre}
For every copula $C$,
\[
|\rho_{11}^{\Lambda}(C)|\le R^*\bigl(\rho_{22}^{\Lambda}(C)\bigr),
\]
and the inequality is sharp for every value of $\rho_{22}^{\Lambda}(C)$.
Equivalently, $\rho_{22}^{\Lambda}(C)\ge(R^*)^{-1}(|\rho_{11}^{\Lambda}(C)|)$
whenever $|\rho_{11}^{\Lambda}(C)|\ge\frac12$, while every value
$\rho_{22}^{\Lambda}\in[-\frac78,1]$ is compatible with
$|\rho_{11}^{\Lambda}|\le\frac12$.
\end{corollary}

\begin{table}[ht]
\centering
\caption{Smallest arachnitude compatible with a given value of Spearman's
rho, $(R^*)^{-1}(|\rho_{11}^{\Lambda}|)$, computed from the arcs of
\cref{thm:main-legendre}.}
\label{tab:min-arachnitude}
\small\setlength{\tabcolsep}{4pt}

\medskip
\begin{tabular}{lcccccccc}
\toprule
$|\rho_{11}^{\Lambda}|$ & $\le0.5$ & $0.55$ & $0.6$ & $0.7$ & $0.744$ & $0.8$ & $0.9$ & $0.95$\\
\midrule
$\min\rho_{22}^{\Lambda}$ & $-0.875$ & $-0.806$ & $-0.723$ & $-0.505$ & $-0.382$ & $-0.197$ & $0.249$ & $0.558$\\
\bottomrule
\end{tabular}
\end{table}

\Cref{cor:sharp-legendre} is a consistency condition for first- and
second-order rank dependence: a pair of estimated or expert-specified values
$(\rho_{22}^{\Lambda},\rho_{11}^{\Lambda})$ outside $\A^{\Lambda}$ cannot be
produced by any copula. For instance, \cref{tab:min-arachnitude} shows that
a strongly concordant pair of risks cannot be strongly anti-cruciform.

\begin{figure}[t]
\centering
\includegraphics[width=\textwidth]{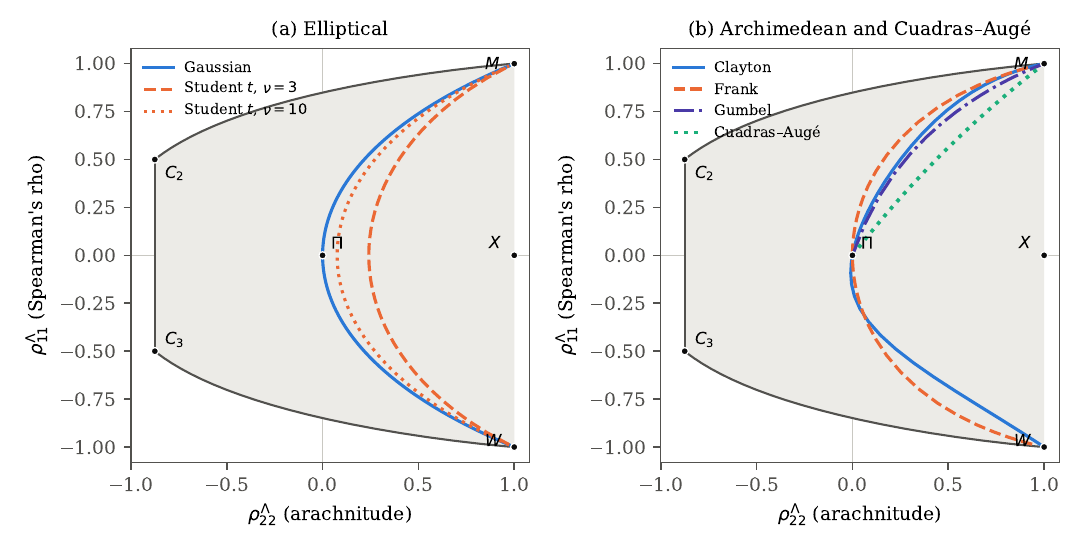}
\caption{The region $\A^{\Lambda}$ of (arachnitude, Spearman's rho) (shaded;
see \cref{thm:main-legendre}) with parametric families. (a) Gaussian and
Student~$t$ copulas, $r\in(-1,1)$. (b) Clayton, Frank and Gumbel copulas,
and Cuadras--Aug\'e copulas. Points: $M$, $W$, independence $\Pi$, the cross
copula $X$ and the minimum-arachnitude copulas $C_2,C_3$.}
\label{fig:families}
\end{figure}

We conclude the section with the values of the two coefficients for some
copulas and classical families, which are represented in
\cref{fig:families}. These values were computed by numerical integration of
$\rhoo(C)$ and of $\rhot(C)$, with $\ARV$ given by
\eqref{eq:arv-lebesgue}, for copulas with a closed form, and of
$\E[\Lambda_j(U)\Lambda_j(V)]$ for elliptical copulas. We start from the
extreme values of arachnitude.

\begin{example}[Extreme values of arachnitude]\label{ex:maximum}
$M$, $W$ and the cross copula $X=(M+W)/2$ have folded copula $M$, hence
$\rhot=1$, while $\rhoo=1,-1,0$, respectively. At the other extreme, the
shuffles of Min $C_2$ and $C_3$ of \cref{subsec:familyI} have folded copula
$W$, hence $\rhot=-\frac78$, while $\rhoo=\frac12$ and $\rhoo=-\frac12$; for
these copulas one rank is extreme exactly when the other is central. The
mixtures $(1-p)\Pi+pX$ and $(1-p)\Pi+p\frac{C_2+C_3}{2}$, $p\in[0,1]$, have
$\rhoo=0$ and $\rhot=p$, respectively $\rhot=-\frac78p$; they fill the
horizontal axis of the region.
\end{example}

Among elliptical copulas, the Gaussian family admits an explicit
description.

\begin{proposition}\label{prop:gaussian}
Let $C_r$ be the Gaussian copula with correlation $r\in[-1,1]$. Then
$\rhoo(C_r)=\frac6\pi\arcsin\frac r2$ and
\[
\rhot(C_r)=\sum_{k\ge1}\beta_{2k}\,r^{2k},\qquad
\beta_{2k}=\frac{\bigl(\E[\Lambda_2(\Phi(Z))\,He_{2k}(Z)]\bigr)^2}{(2k)!}\ge0,
\]
with $\sum_{k\ge1}\beta_{2k}=1$, where $Z$ is standard normal, $\Phi$ its distribution function and $He_m$ the
probabilists' Hermite polynomials. In particular
$\beta_2=\frac{15}{2\pi^2}\approx0.760$, so
$\rhot(C_r)=\frac{15}{2\pi^2}r^2+O(r^4)$, and $\rhot(C_r)$ is even,
non-negative and increasing in $|r|$.
\end{proposition}

\begin{proof}
The formula for $\rhoo$ is classical. Write $f=\Lambda_2\circ\Phi$, so
$\E f(Z)=0$ and $\E f(Z)^2=1$. Its Hermite expansion
$f=\sum_m c_m He_m/m!$, $c_m=\E[f(Z)He_m(Z)]$, has only even terms because
$f$ is even. By Mehler's formula (see, e.g., \cite{Lancaster1958}), for
$(Z_1,Z_2)$ standard bivariate normal with correlation $r$,
$\E[f(Z_1)f(Z_2)]=\sum_m c_m^2r^m/m!$, and $\sum_m c_m^2/m!=\E f(Z)^2=1$.
For $m=2$, Stein's identity gives $c_2=\E[f(Z)(Z^2-1)]=\E[f''(Z)]$, where,
with $\phi=\Phi'$,
\[
f''(z)=\sqrt5\bigl(12\phi(z)^2-12z\Phi(z)\phi(z)+6z\phi(z)\bigr).
\]
Now $\E[Z\phi(Z)]=0$ by symmetry, $\E[Z\Phi(Z)\phi(Z)]=\frac12\E[\phi(Z)^2]$
by Stein's identity $\E[Zk(Z)]=\E[k'(Z)]$ applied to $k=\Phi\phi$, and
$\E[\phi(Z)^2]=\frac1{2\pi\sqrt3}$. Hence
$c_2=6\sqrt5\,\E[\phi(Z)^2]=\sqrt{15}/\pi$ and
$\beta_2=c_2^2/2=15/(2\pi^2)$.
\end{proof}

Combining the two leading terms, Gaussian copulas satisfy
$\rhot\approx\frac56(\rhoo)^2$ for weak dependence. A sample with small
Spearman's rho and clearly positive or negative arachnitude is therefore not
Gaussian. Numerically, $\rhot(C_{0.5})\approx0.202$ and
$\rhot(C_{0.9})\approx0.763$.

Student~$t$ copulas behave differently, since they exhibit joint extremes
even without correlation.

\begin{proposition}\label{prop:t}
Let $C_{r,\nu}$ be the Student~$t$ copula with correlation parameter $r$ and
$\nu>0$ degrees of freedom. Then $\rhot(C_{0,\nu})>0$, although
$\rhoo(C_{0,\nu})=0$.
\end{proposition}

\begin{proof}
Write $(X,Y)=\sqrt{G}\,(Z_1,Z_2)$, with $Z_1,Z_2$ independent standard normal
and $G=\nu/\chi^2_\nu$ independent of them. By symmetry of the $t$
distribution, $(2U-1)^2=k(|X|)$ and $(2V-1)^2=k(|Y|)$ with $k$ strictly
increasing. Conditionally on $G$, $|X|$ and $|Y|$ are independent, so
\[
\Cov\bigl(k(|X|),k(|Y|)\bigr)=\Var\bigl(m(G)\bigr),\qquad
m(g)=\E\bigl[k(\sqrt g\,|Z_1|)\bigr],
\]
and $m$ is strictly increasing. Since $G$ is not degenerate, the variance is
positive. Finally $\rhoo=0$ by the symmetry $(X,Y)\mapsto(X,-Y)$.
\end{proof}

The mixing variable $G$ inflates and deflates both coordinates together,
which produces the cross-shaped clustering in all four corners that
arachnitude detects. Numerically $\rhot(C_{0,3})\approx0.241$ and
$\rhot(C_{0,10})\approx0.076$. For fixed $r$ the $t$ curves lie to the right
of the Gaussian one (\cref{fig:families}a): at the same Spearman's rho, a
$t$ copula has larger arachnitude.

We now turn to Archimedean copulas.

For the Clayton, Frank and Gumbel families (\cref{fig:families}b) the curves
start at $\Pi$ and end at $M$ (and at $W$ for Clayton and Frank). All three
stay in, or extremely close to, the half-plane $\rhot\ge0$; the smallest value
we found is $\rhot\approx-0.009$ for Clayton copulas with negative parameter.
At equal Spearman's rho, Frank copulas, which have no tail dependence, have
the smallest arachnitude, while Clayton and Gumbel copulas, with lower and
upper tail dependence respectively, have larger values.

A class of copulas for which the two coefficients can be computed in closed
form is that of lower semilinear copulas
\cite{DuranteKolesarovaMesiarSempi2008}. These copulas are linear on each
segment joining a point of the main diagonal to the lower or to the left
boundary of the unit square; equivalently, they can be written as
\[
C_f(u,v)=\min\{u,v\}\,f(\max\{u,v\}),
\]
where $f:(0,1]\to\R$ is non-decreasing, $t\mapsto f(t)/t$ is
non-increasing and $f(1)=1$ (here $f(t)=\delta(t)/t$, $\delta$ being the
diagonal section). Since $t\le f(t)\le1$, every lower semilinear copula
satisfies $\Pi\le C_f\le M$.

\begin{proposition}\label{prop:lsc}
Let $C_f$ be a lower semilinear copula and let $\mu_f$ be the non-negative
Borel measure on $(0,1]$ such that $f(t)/t-1=\mu_f((t,1])$ for all but
countably many $t$. Then
\begin{align*}
\rhoo(C_f)&=12\int_0^1t^2f(t)\,dt-3=3\int_{(0,1]}s^4\,d\mu_f(s),\\
\rhot(C_f)&=45\int_{(0,1]}s^4\Bigl(1-\frac{4s}3\Bigr)^2d\mu_f(s).
\end{align*}
Consequently, $0\le\rhot(C_f)\le15\,\rhoo(C_f)$, and $\rhot(C_f)=0$ or
$\rhoo(C_f)=0$ only if $C_f=\Pi$.
\end{proposition}

\begin{proof}
Since $C_f$ is symmetric, $\int C_f=2\int_0^1f(y)\int_0^yx\,dx\,dy
=\int_0^1y^2f(y)\,dy$, and, by \eqref{eq:arv-lebesgue},
$\ARV(C_f)=8\int_0^1f(y)(1-2y)\int_0^yx(1-2x)\,dx\,dy
=\frac43\int_0^1y^2(1-2y)(3-4y)f(y)\,dy$. Write $f(t)=t+t\,g(t)$ with
$g(t)=\mu_f((t,1])$. The term $f(t)=t$ corresponds to $\Pi$, for which both
coefficients vanish. For the second term, Fubini's theorem gives
$\int_0^1w(t)\,t\,g(t)\,dt=\int_{(0,1]}\int_0^sw(t)\,t\,dt\,d\mu_f(s)$, and
\[
\int_0^s12t^3\,dt=3s^4,\qquad
\int_0^s\tfrac43t^3(1-2t)(3-4t)\,dt=s^4\Bigl(1-\frac{4s}3\Bigr)^2 .
\]
With \cref{thm:arv-arachnitude}, this gives the two formulas. The
inequalities follow from $0\le(1-\frac{4s}3)^2\le1$ on $(0,1]$. If
$\rhoo(C_f)=0$, then $\mu_f=0$; if $\rhot(C_f)=0$, then $\mu_f$ is
concentrated at $s=\frac34$, which is compatible with the monotonicity of
$f$ only if $\mu_f=0$. In both cases $f(t)=t$ and $C_f=\Pi$.
\end{proof}

Thus lower semilinear copulas, which are positively quadrant dependent,
always have non-negative arachnitude: they cannot be anti-cruciform. The most
prominent example is the Cuadras--Aug\'e family \cite{CuadrasAuge1981},
which is also a family of extreme-value copulas.

\begin{example}[Cuadras--Aug\'e copulas]\label{ex:CA}
For $\theta\in[0,1]$, the Cuadras--Aug\'e copula
$C_\theta(u,v)=\min\{u,v\}^\theta(uv)^{1-\theta}$ is lower semilinear with
$f(t)=t^{1-\theta}$. By \cref{prop:lsc}, or by direct integration,
\begin{gather*}
\rhoo(C_\theta)=\frac{3\theta}{4-\theta},
\qquad
\rhot(C_\theta)=\frac{5\theta(\theta^2-3\theta+14)}{(4-\theta)(5-\theta)(6-\theta)},\\
\ARV(C_\theta)=\frac19+\frac{\theta(\theta^2-3\theta+14)}{9(4-\theta)(5-\theta)(6-\theta)} .
\end{gather*}
The curve joins $\Pi$ ($\theta=0$) to $M$ ($\theta=1$). For small $\theta$,
$\rhot(C_\theta)=\frac{7}{12}\theta+O(\theta^2)$ and
$\rhoo(C_\theta)=\frac34\theta+O(\theta^2)$, so that
$\rhot\approx\frac79\rhoo$: arachnitude grows linearly with Spearman's rho,
whereas for Gaussian copulas it grows quadratically
(\cref{prop:gaussian}). This reflects the singular component of
$C_\theta$ on the main diagonal, which places mass at points where both
ranks are simultaneously extreme (\cref{fig:families}b).
\end{example}

\section{Concluding remarks}
\label{sec:conclusions}

We have determined the exact $(\ARV,\rho)$ region of bivariate copulas.
Equivalently, the sharp inequality
\[
|\rho(C)|\le r^*\bigl(\ARV(C)\bigr),
\qquad C\in\mathcal C,
\]
holds, with equality on the upper boundary for
$C^{\mathrm I}_a$ ($0\le a\le\frac25$),
$C^{\mathrm{II}}_{b,c(b)}$ ($0\le b\le\frac15$) and
$C^{\mathrm{III}}_c$ ($c^*\le c\le1$), and on the lower boundary for
their flipped versions. Since $r^*$ is continuous and strictly
increasing from $\frac12$ to $1$, the same result can be read as a sharp
lower bound for $\ARV$ in terms of Spearman's rho:
\[
\ARV(C)\ge(r^*)^{-1}\bigl(|\rho(C)|\bigr)
\qquad\text{whenever } |\rho(C)|\ge\tfrac12,
\]
while for $|\rho(C)|\le\frac12$ only the global bound
$\ARV(C)\ge\frac{11}{120}$ is available, and it is attained. The upper boundary is
continuously differentiable away from its left endpoint
(Corollary~\ref{cor:C1}) and has a vertical tangent there
(Corollary~\ref{cor:vertical-tangent}).

By \cref{sec:arachnitude}, on copulas the average rectangular volume is an
affine function of arachnitude, $\rhot=45\,\ARV-5$. An index introduced
for aggregation functions through rectangular volumes thus turns out to be
a second-order rank correlation, the natural companion of Spearman's rho,
and \cref{thm:main} is equivalently a sharp compatibility condition between
first- and second-order rank dependence (\cref{thm:main-legendre}).

Two ingredients of the proof may be useful beyond the present problem.
First, the folding transformation turns $\ARV$ into a moment of the
folded copula and gives a sharp upper bound for $|\rho|$ in terms of
another moment; the same reduction applies to every pair of functionals
that, after folding, become moments of the form $\E[g(S)h(T)]$.
Second, the supporting costs $c_\kappa$ are neither supermodular nor
submodular for $0<\kappa<2$, so the optimal couplings are not monotone
rearrangements. They were obtained by constructing explicit Kantorovich
potentials, and the algebraic parameter constraint $P(b,c)=0$ appeared
as the continuity condition of the potential.

\section*{Declaration of generative AI and AI-assisted technologies in the manuscript preparation process}

During the preparation of this work the authors used Claude (Anthropic) in order to explore the problem numerically, to propose and check constructions and proofs (in particular, related to Families~II and~III with the associated Kantorovich potentials), to produce
the figures, and to improve the organization and the language of the manuscript. After using this tool, the authors reviewed and edited the content as needed and take full responsibility for the content of the published article.

\section*{Declaration of competing interest}
The authors declare that they have no known competing financial interests or
personal relationships that could have appeared to influence the work
reported in this paper.

\section*{Data availability}
No data were used for the research described in the article. The Python
code that reproduces all figures is provided as supplementary material.

\section*{Acknowledgments}
The authors are members of the group GNAMPA of INdAM, which is
acknowledged.

\appendix

\renewcommand{\thetheorem}{\Alph{section}.\arabic{theorem}}
\renewcommand{\theequation}{\Alph{section}.\arabic{equation}}
\renewcommand{\thefigure}{\Alph{section}.\arabic{figure}}
\renewcommand{\thetable}{\Alph{section}.\arabic{table}}
\renewcommand{\thesubsection}{\Alph{section}.\arabic{subsection}}
\renewcommand{\thesubsubsection}{\thesubsection.\arabic{subsubsection}}

\section{Proof of the main result}
\label{app:proofs}

This section contains the proof of \cref{thm:main} and of
\cref{cor:C1,cor:vertical-tangent}. By the folding transformation of
\cref{prop:folding}, the upper boundary is described by the maximization of
$\E_D[ST]$ over copulas $D$ with a prescribed value of $\E_D[S^2T^2]$;
supporting lines lead to a family of optimal transport problems. For each
family we exhibit a Kantorovich potential, which proves optimality. For
Families~II and~III the potential is continuous exactly when $P(b,c)=0$,
which selects the parameter curve $c=c(b)$.

\subsection{\texorpdfstring{Supporting costs and Kantorovich duality}{Supporting costs and Kantorovich duality}}
\label{subsec:transport}

For $\kappa>0$, define
\begin{equation}
	h_\kappa(z)
	=
	z-\frac{z^2}{2\kappa},
	\qquad
	c_\kappa(s,t)
	=
	h_\kappa(st)
	=
	st-\frac{s^2t^2}{2\kappa}.
	\label{eq:transport-cost}
\end{equation}

Write $x=\E_D[S^2T^2]$ and $y=\E_D[ST]$. The set of pairs $(x,y)$, $D\in\mathcal C$, is compact and convex, and by
\eqref{eq:rstar-folded} its upper boundary is mapped onto the graph of $r^*$ by the affine map $(x,y)\mapsto(\frac1{12}+\frac x4,3y)$.
A supporting line of positive finite slope $1/(2\kappa)$ in the
$(x,y)$-plane is obtained by maximizing $y-x/(2\kappa)$, that is, by
solving
\begin{equation}
	\sup_{D\in\mathcal C}\E_D[c_\kappa(S,T)].
	\label{eq:transport-problem}
\end{equation}
In the $(\ARV,\rho)$-plane the corresponding slope is $6/\kappa$.
We shall construct optimizers for suitable values of $\kappa$ and
then show that their moment pairs cover the upper boundary.
See, for example, \cite{RachevRuschendorf1998,Villani2009} for the
general theory of optimal transport.

For our purposes, the following elementary form of Kantorovich
duality is sufficient.

\begin{lemma}[Kantorovich duality: a sufficient optimality condition]\label{lem:kantorovich-sufficient}
	Let $D\in\mathcal C$ and $\kappa>0$. Suppose that there exists a bounded
	Borel function $\phi:[0,1]\to\mathbb R$ such that
	\[
	G(s,t)
	:=
	\phi(s)+\phi(t)-c_\kappa(s,t)
	\ge0
	\qquad
	\text{for all }(s,t)\in[0,1]^2,
	\]
	and
	\[
	G=0
	\qquad
	D\text{-almost surely}.
	\]
	Then $D$ solves \eqref{eq:transport-problem}. In particular, if
	\[
	\alpha_D
	=
	\frac1{12}
	+
	\frac14\,\E_D[S^2T^2],
	\]
	then
	\[
	r^*(\alpha_D)
	=
	3\,\E_D[ST].
	\]
\end{lemma}

\begin{proof}
	For every $D'\in\mathcal C$,
	\[
	\E_{D'}[c_\kappa(S,T)]
	\le
	\E_{D'}[\phi(S)+\phi(T)]
	=
	2\int_0^1\phi(u)\,du.
	\]
	Since equality holds for $D$, this proves the optimality of $D$.
	If $D'$ has the same value of $\E[S^2T^2]$ as $D$, the inequality
	\[
	\E_{D'}[c_\kappa]
	\le
	\E_D[c_\kappa]
	\]
	reduces to
	\[
	\E_{D'}[ST]
	\le
	\E_D[ST].
	\]
	The conclusion follows from \eqref{eq:rstar-folded}.
\end{proof}

\paragraph{Terminology.} Lemma~\ref{lem:kantorovich-sufficient} is the elementary
half of Kantorovich duality; see \cite[Ch.~5]{Villani2009} and
\cite[Ch.~1]{Santambrogio2015}. The dual problem associated with
\eqref{eq:transport-problem} is
\[
\inf\Bigl\{\int_0^1\phi(u)\,du+\int_0^1\psi(u)\,du:\
\phi(s)+\psi(t)\ge c_\kappa(s,t)\text{ for all }(s,t)\in[0,1]^2\Bigr\},
\]
and, since $c_\kappa$ is symmetric and both marginals are uniform, it
suffices to consider symmetric pairs $\psi=\phi$. A function $\phi$ as
in Lemma~\ref{lem:kantorovich-sufficient} is called a
\emph{Kantorovich potential} (for $D$ and the cost $c_\kappa$): the
pair $(\phi,\phi)$ solves the dual problem and the primal and dual
values coincide. We call $G$ the \emph{duality slack}; the requirement
$G=0$ $D$-almost surely is the \emph{complementary slackness} condition.
The set
\[
\Gamma_\phi=\{(s,t)\in[0,1]^2:\ G(s,t)=0\}
\]
is the \emph{contact set} of $\phi$, and its section
\[
\partial^{c}\phi(t)=\{s\in[0,1]:\ (s,t)\in\Gamma_\phi\}
\]
is the $c$-superdifferential of $\phi$ at $t$ (the convention of
\cite[Ch.~5]{Villani2009}, adapted to a maximization problem). If
$\phi$ is continuous, $\Gamma_\phi$ is closed and complementary
slackness is equivalent to $\operatorname{supp}(D)\subseteq\Gamma_\phi$.

For a copula $D$ and $t\in[0,1]$ we write
\[
\Gamma_D(t)=\{p\in[0,1]:\ (p,t)\in\operatorname{supp}(D)\}
\]
for the set of points \emph{coupled with} $t$ by $D$. Under the
assumptions of Lemma~\ref{lem:kantorovich-sufficient} with $\phi$
continuous, $\Gamma_D(t)\subseteq\partial^{c}\phi(t)$.

If $(p,t)\in\Gamma_\phi$ with $0<t<1$ and $\phi$ is differentiable at
$t$, then $t$ minimizes the non-negative function $G(p,\cdot)$, hence
$\partial_tG(p,t)=0$, that is,
\[
\phi'(t)=\partial_t c_\kappa(p,t)=p\,h_\kappa'(pt).
\]
This first-order condition is the one-dimensional counterpart of the
relation between the gradient of a Kantorovich potential and the
optimal transport map; cf.\ \cite{GangboMcCann1996}.

We shall repeatedly use the following factorization of
$\partial_tG$. For $t>0$, let $p\in\Gamma_D(t)$, assume that $\phi$ is differentiable at $t$ and
\[
\phi'(t)
=
p\,h_\kappa'(pt),
\]
then
\begin{align}
	\frac{\partial G}{\partial t}(s,t)
	&=
	p\,h_\kappa'(pt)
	-
	s\,h_\kappa'(st)
	\nonumber\\
	&=
	(p-s)
	\left(
	1-\frac{t(p+s)}{\kappa}
	\right)
	\nonumber\\
	&=
	\frac{t}{\kappa}
	(s-p)(s-q),
	\qquad
	q=\frac{\kappa}{t}-p.
	\label{eq:slack-factorization}
\end{align}
Thus the sign of $\partial_tG$ is determined by the position of $s$
relative to the two points $p$ and $q$ (the second root $q$ need not be
coupled with $t$).

\subsection{Extremality of Family I}
\label{sec:familyI-optimality}

We first prove that Family~I generates the rightmost arc
$\Gamma^{\mathrm I}$ of the upper boundary. By
\eqref{eq:fold-familyI}, the folded copula of $C^{\mathrm I}_a$ is
$D_b$ with
\begin{equation}
	b=1-2a.
	\label{eq:b-a}
\end{equation}
Thus the range $0\le a\le2/5$ corresponds to
\[
\frac15\le b\le1.
\]

\subsubsection{Optimality for \texorpdfstring{$b\ge1/5$}{b >= 1/5}}
\label{subsec:familyI-optimality}

For $b\in[1/5,1]$, set
\begin{equation}
	\kappa_b
	=
	\frac{(1+b)(1+5b)}{6}.
	\label{eq:kappa-b}
\end{equation}
The corresponding transport cost is
\[
c_{\kappa_b}(s,t)
=
h_{\kappa_b}(st)
=
st-\frac{s^2t^2}{2\kappa_b}.
\]

\begin{theorem}\label{thm:familyI-optimal}
	For every $a\in[0,2/5]$,
	\[
	\left(
	\ARV(C^{\mathrm I}_a),
	\rho(C^{\mathrm I}_a)
	\right)
	\]
	lies on the upper boundary of $\A$. Equivalently,
	$\Gamma^{\mathrm I}$ is a boundary arc.
\end{theorem}

\begin{proof}
	The case $a=0$ is immediate, since $C^{\mathrm I}_0=M$.
	Assume therefore $a>0$ and let $b=1-2a\in[1/5,1)$.
	
	We construct a Kantorovich potential for $D_b$. Define
	$\phi_b:[0,1]\to\mathbb R$ by
	\begin{equation}
		\phi_b(s)
		=
		\frac12 h_{\kappa_b}(s^2),
		\qquad
		0\le s\le b,
		\label{eq:phi-b-left}
	\end{equation}
	and
	\begin{equation}
		\phi_b(s)
		=
		\phi_b(b)
		+
		\int_b^s
		(1+b-y)\,
		h_{\kappa_b}'\!\bigl(y(1+b-y)\bigr)\,dy,
		\qquad
		b\le s\le1.
		\label{eq:phi-b-right}
	\end{equation}
	Set
	\[
	G_b(s,t)
	=
	\phi_b(s)+\phi_b(t)-c_{\kappa_b}(s,t).
	\]
	On the diagonal piece, $G_b=0$ by definition. On the anti-diagonal,
	the derivative of $s\mapsto G_b(s,1+b-s)$ is zero.
	Evaluating the integral in \eqref{eq:phi-b-right} at $s=1$ gives
	\[
	\phi_b(b)+\phi_b(1)-h_\kappa(b)
	=
	\frac{(1-b)^2}{12\kappa}
	\bigl(6\kappa-(1+b)(1+5b)\bigr).
	\]
	For $\kappa=\kappa_b$, the right-hand side vanishes by
	\eqref{eq:kappa-b}. Hence
	\[
	\phi_b(b)+\phi_b(1)=h_{\kappa_b}(b).
	\]
	Hence the constant on the anti-diagonal is also zero. Thus $G_b$ vanishes on
	\[
	\operatorname{supp}(D_b)
	=
	\{(s,s):0\le s\le b\}
	\cup
	\{(s,t)\in[b,1]^2:s+t=1+b\}.
	\]
	It remains to prove that $G_b\ge0$ on $[0,1]^2$.
	
	\medskip
	\noindent\emph{The block $[0,b]^2$.}
	A direct expansion gives
	\begin{equation}
		G_b(s,t)
		=
		\frac{(s-t)^2}{2}
		\left(
		1-\frac{(s+t)^2}{2\kappa_b}
		\right).
		\label{eq:G-diagonal-block}
	\end{equation}
	Since
	\[
	(s+t)^2\le4b^2\le2\kappa_b,
	\]
	we obtain $G_b(s,t)\ge0$. Indeed,
	\[
	4b^2\le2\kappa_b
	\]
	is equivalent to
	\[
	12b^2\le(1+b)(1+5b),
	\]
	which holds for $0\le b\le1$.
	
	\medskip
	\noindent\emph{The block $[b,1]^2$.}
	For $(s,t)\in[b,1]^2$ one obtains
	\begin{equation}
		G_b(s,t)
		=
		-\frac{(s+t-1-b)^2}
		{2(1+b)(1+5b)}
		\,Q_b(s,t),
		\label{eq:G-antidiagonal-block}
	\end{equation}
	where
	\begin{equation}
		Q_b(s,t)
		=
		3(s-t)^2
		-2(1+b)(s+t)
		+4b(1+b).
		\label{eq:Qb}
	\end{equation}
	Write
	\[
	x=s+t,
	\qquad
	y=s-t.
	\]
	If $x\le1+b$, then $|y|\le x-2b$, and therefore
	\begin{align*}
		Q_b(s,t)
		&\le
		3(x-2b)^2
		-2(1+b)(x-2b)\\
		&=
		(x-2b)
		\bigl[3(x-2b)-2(1+b)\bigr]
		\le0,
	\end{align*}
	because
	\[
	0\le x-2b\le1-b
	\]
	and
	\[
	3(1-b)\le2(1+b)
	\quad\Longleftrightarrow\quad
	b\ge\frac15.
	\]
	If $x\ge1+b$, then $|y|\le2-x$. Setting $z=2-x$, so that
	$0\le z\le1-b$, we obtain
	\[
	Q_b(s,t)
	\le
	3z^2+2(1+b)z-4(1-b^2).
	\]
	The right-hand side is increasing in $z\ge0$, and at $z=1-b$ it equals
	\[
	(1-b)(1-5b)\le0.
	\]
	Hence $Q_b\le0$ also in this case, and
	\eqref{eq:G-antidiagonal-block} yields $G_b\ge0$ on $[b,1]^2$.
	
	\medskip
	\noindent\emph{The mixed blocks.}
	By symmetry it is enough to consider $0\le s\le b\le t\le1$.
	For fixed $t$, let
	\[
	F(s)=G_b(s,t).
	\]
	From \eqref{eq:phi-b-left},
	\[
	F''(s)
	=
	1-\frac{3s^2}{\kappa_b}
	+\frac{t^2}{\kappa_b}
	\ge
	1-\frac{2b^2}{\kappa_b}
	\ge0.
	\]
	Moreover,
	\begin{align*}
		F'(b^-)
		&=
		b\,h_{\kappa_b}'(b^2)
		-
		t\,h_{\kappa_b}'(bt)\\
		&=
		(b-t)
		\left(
		1-\frac{b(b+t)}{\kappa_b}
		\right)
		\le0,
	\end{align*}
	because
	\[
	\kappa_b\ge b(1+b).
	\]
	Since $F$ is convex, $F'$ is non-decreasing on $[0,b]$, hence $F'(s)\le F'(b^-)\le0$ for $s\le b$.
	Therefore
	\[
	G_b(s,t)=F(s)\ge F(b)=G_b(b,t)\ge0,
	\]
	where the last inequality follows from the already treated block
	$[b,1]^2$.
	
	We have thus proved
	\[
	G_b\ge0
	\qquad\text{on }[0,1]^2,
	\]
	with equality $D_b$-almost surely. Hence $\phi_b$ is a Kantorovich
	potential for $D_b$. Lemma~\ref{lem:kantorovich-sufficient} shows that $D_b$
	is optimal for the corresponding transport problem. Since
	\[
	\bigl(C^{\mathrm I}_a\bigr)^{\mathsf F}=D_b,
	\]
	Lemma~\ref{lem:kantorovich-sufficient}, \eqref{eq:fold-familyI} and $\rho(\unfold{D_b})=3\E_{D_b}[ST]$ give
	\[
	r^*\!\left(\ARV(C^{\mathrm I}_a)\right)
	=
	\rho(C^{\mathrm I}_a),
	\]
	as claimed.
\end{proof}

\begin{remark}\label{rem:kappa-familyI}
	The value of $\kappa_b$ is determined by the tangent to the curve
	generated by Family~I. More precisely, for $0\le b<1$,
	\[
	\frac{1}{2\kappa_b}
	=
	\frac{3}{(1+b)(1+5b)}
	=
	\frac{\dfrac{d}{db}\E_{D_b}[ST]}
	{\dfrac{d}{db}\E_{D_b}[S^2T^2]}.
	\]
	Thus the cost used in the proof is the supporting cost associated
	with the tangent to $\Gamma^{\mathrm I}$.
\end{remark}

\subsubsection{Sharpness of the threshold}
\label{subsec:familyI-threshold}

The restriction $a\le2/5$ in
Theorem~\ref{thm:familyI-optimal} is sharp.

\begin{proposition}\label{prop:familyI-not-optimal}
	For
	\[
	\frac25<a<\frac12,
	\]
	the point
	\[
	\left(
	\ARV(C^{\mathrm I}_a),
	\rho(C^{\mathrm I}_a)
	\right)
	\]
	does not lie on the upper boundary of $\A$.
\end{proposition}

\begin{proof}
	Let
	\[
	b=1-2a,
	\]
	so that $0<b<1/5$, and let $\kappa_b$ be given by
	\eqref{eq:kappa-b}. The mixed derivative of the corresponding cost is
	\[
	\frac{\partial^2 c_{\kappa_b}}
	{\partial s\,\partial t}(s,t)
	=
	1-\frac{2st}{\kappa_b}.
	\]
	At the corner $(b,1)$ of the anti-diagonal part of
	$\operatorname{supp}(D_b)$,
	\[
	\kappa_b-2b
	=
	\frac{(1-b)(1-5b)}{6}
	>0.
	\]
	Hence
	\[
	\frac{\partial^2 c_{\kappa_b}}
	{\partial s\,\partial t}(b,1)>0,
	\]
	and the cost is strictly supermodular in a neighborhood of $(b,1)$.
	
	We make the local improvement explicit. Choose two disjoint
	intervals $I$ and $I+d$ of equal positive length in a sufficiently
	small right neighborhood of $b$. On these intervals the original
	map is $T(s)=1+b-s$. Replace it by $T(s+d)$ for $s\in I$ and
	by $T(s-d)$ for $s\in I+d$, leaving it unchanged elsewhere.
	This exchanges the two image intervals, so both uniform marginals
	are preserved.
	
	For $s\in I$, put $x_1=s$, $x_2=s+d$,
	$y_1=T(s+d)$ and $y_2=T(s)$. The new matching pairs
	$(x_1,y_1)$ and $(x_2,y_2)$, whereas the old matching pairs
	$(x_1,y_2)$ and $(x_2,y_1)$. Its gain is
	\begin{align*}
		&c_{\kappa_b}(x_1,y_1)+c_{\kappa_b}(x_2,y_2)
		-c_{\kappa_b}(x_1,y_2)-c_{\kappa_b}(x_2,y_1)\\
		&\hspace{12mm}=
		\int_{x_1}^{x_2}\int_{y_1}^{y_2}
		\frac{\partial^2 c_{\kappa_b}}{\partial x\,\partial y}(x,y)
		\,dy\,dx>0.
	\end{align*}
	The intervals can be chosen so that every such rectangle lies in
	the neighborhood where the mixed derivative is positive.
	Integrating over $I$ proves that $D_b$ does not maximize
	$\E_D[c_{\kappa_b}(S,T)]$.
	
	It remains to connect this fact with the boundary problem. In the moment
	plane, put
	\[
	x(b)=\E_{D_b}[S^2T^2],
	\qquad
	y(b)=\E_{D_b}[ST].
	\]
	For $0<b<1/5$ the parameter value is interior and
	\[
	x'(b)
	=
	\frac{(1-b)^2(1+b)(1+5b)}{6}
	>0.
	\]
	If $(x(b),y(b))$ were on the upper boundary of the convex feasible
	moment region, any supporting line there would have to be tangent to
	the differentiable feasible curve $b\mapsto(x(b),y(b))$. Its slope
	would therefore be
	\[
	\frac{y'(b)}{x'(b)}
	=
	\frac{1}{2\kappa_b}
	\]
	by Remark~\ref{rem:kappa-familyI}. Thus $D_b$ would maximize
	\[
	y-\frac{1}{2\kappa_b}x
	=
	\E_D[c_{\kappa_b}(S,T)],
	\]
	contradicting the preceding local improvement. Consequently
	\[
	\left(
	\ARV(C^{\mathrm I}_a),
	\rho(C^{\mathrm I}_a)
	\right)
	\]
	cannot lie on the upper boundary.
\end{proof}

\begin{remark}[An explicit counterexample]
	For $a\in(\frac25,\frac12)$ and $0<\delta\le\frac{5a-2}6$, let $C_{a,\delta}$ be the shuffle of Min obtained from $C^{\mathrm I}_a$ by replacing, in
	$[1-a,1]^2$, the anti-diagonal segment by the three segments $v=u+a-\delta$ ($u\in[1-a,1-a+\delta]$), $v=2-a-u$ ($u\in[1-a+\delta,1-\delta]$),
	$v=u-a+\delta$ ($u\in[1-\delta,1]$), and symmetrically in $[0,a]^2$. Then $\ARV(C_{a,\delta})=\ARV(C^{\mathrm I}_a)+\frac83\delta^3(1-2a+\delta)(1-\delta)$,
	$\rho(C_{a,\delta})=\rho(C^{\mathrm I}_a)+8\delta^3$, and the point $(\ARV(C_{a,\delta}),\rho(C_{a,\delta}))$ lies strictly above the Family~I curve: its vertical distance from the curve is at least
	$8\delta^3(a(5a-2-6\delta)+3\delta^2)/((1-a)(3-5a))>0$ (use the concavity of the Family~I curve and its slope $9/((1-a)(3-5a))$). For $a=\frac9{20}$, $\delta=\frac1{24}$: $\ARV=\frac{901447783}{9331200000}$,
	$\rho=\frac{137393}{216000}\approx0.6360787$, while the Family~I curve gives $\approx0.6360708$ at the same abscissa. This makes
	Proposition~\ref{prop:familyI-not-optimal} concrete.
\end{remark}

\begin{remark}
	The endpoint $a=1/2$ is exceptional. Although the Family~I curve is
	not extremal for $2/5<a<1/2$, one has
	\[
	C^{\mathrm I}_{1/2}=C_2=C^{\mathrm{III}}_1,
	\]
	and therefore
	\[
	\left(
	\ARV(C^{\mathrm I}_{1/2}),
	\rho(C^{\mathrm I}_{1/2})
	\right)
	=
	\left(
	\frac{11}{120},
	\frac12
	\right),
	\]
	which is the left endpoint of the upper boundary. In the decomposition
	of Theorem~\ref{thm:main}, this point belongs to the arc generated by
	Family~III.
\end{remark}

\subsection{Extremality of Families II and III}
\label{sec:familiesIIIII-optimality}

We now prove the extremality of the two new families introduced in
Section~\ref{sec:families}. After folding, both are described by the
same family $D_{b,c}$: Family~II corresponds to $b>0$, whereas
Family~III is the limiting case $b=0$.

The proof has three ingredients. First we record the geometry of the
support of $D_{b,c}$. We then construct a candidate Kantorovich
potential and show that the polynomial $P(b,c)$ introduced in
Section~\ref{sec:main} is precisely the condition that makes this
potential continuous. Finally, a monotonicity argument reduces the
global inequality for the Kantorovich slack to its values on the
support and at one junction point.

\subsubsection{Geometry of the folded family}
\label{subsec:Dbc-geometry}

Throughout this section, let
\begin{equation}
	\Delta
	=
	\left\{
	(b,c):
	0\le b\le\frac15,\quad
	\frac12+b\le c\le1
	\right\}
	\setminus\{(0,1)\},
	\label{eq:Delta}
\end{equation}
and set
\begin{equation}
	\kappa=c(1+2b-c),
	\qquad
	b'=1+b-c.
	\label{eq:kappa-bprime}
\end{equation}
Then $\kappa>0$ on $\Delta$, and
\begin{align}
	\kappa-b(1+b)
	&=(c-b)(1+b-c)\ge0,
	\label{eq:F1}\\
	\kappa-2b
	&=(1-c)(c-2b)\ge0,
	\label{eq:F2}\\
	c^2-\kappa
	&=c(2c-1-2b)\ge0,
	\label{eq:F3}\\
	u(1-u)
	&\le c(1-c)\le\kappa,
	\qquad c\le u\le1.
	\label{eq:F4}
\end{align}

Recall that
\[
\xi_1(t)
=
\frac12\left(\frac{\kappa}{t}+t-1\right),
\qquad
\xi_2(t)
=
\frac12\left(\frac{\kappa}{t}-t+1\right),
\qquad c\le t\le1.
\]
Their endpoint values are
\begin{equation}
	\xi_1(c)=b,
	\qquad
	\xi_2(c)=b',
	\qquad
	\xi_1(1)=\xi_2(1)=\frac{\kappa}{2}.
	\label{eq:xi-endpoints}
\end{equation}
Moreover, $\xi_1$ is increasing (since $\xi_1'(t)=\frac12(1-\kappa/t^2)\ge0$ by \eqref{eq:F3}), $\xi_2$ is decreasing, and
\[
0\le b\le\frac{\kappa}{2}\le b'\le c.
\]

Let
\[
\tau:[b,b']\longrightarrow[c,1]
\]
denote the inverse of $\xi_1$ on $[b,\kappa/2]$ and of $\xi_2$ on
$[\kappa/2,b']$. Explicitly,
\[
\tau(s)
=
\frac12
\left(
1+2s+\sqrt{(1+2s)^2-4\kappa}
\right),
\qquad
b\le s\le\frac{\kappa}{2},
\]
and
\[
\tau(s)
=
\frac12
\left(
1-2s+\sqrt{(1-2s)^2+4\kappa}
\right),
\qquad
\frac{\kappa}{2}\le s\le b'.
\]

Along the two curved branches,
\begin{equation}
	t\xi_1(t)
	=
	\frac12\bigl(\kappa-t(1-t)\bigr)
	\le\frac{\kappa}{2},
	\qquad
	t\xi_2(t)
	=
	\frac12\bigl(\kappa+t(1-t)\bigr)
	\ge\frac{\kappa}{2}.
	\label{eq:branch-products}
\end{equation}

As in Section~\ref{subsec:transport}, we write $\Gamma(t)=\Gamma_{D_{b,c}}(t)$ for the set of points coupled with $t$ (the relation is symmetric because $D_{b,c}$ is). Away from the junction
points, these sets are
\begin{center}
	\begin{tabular}{c|c}
		range of $t$ & $\Gamma(t)$\\
		\hline
		$0<t<b$ & $t$\\
		$b<t<b'$ & $\tau(t)$\\
		$b'<t<c$ & $1+b-t$\\
		$c<t<1$ & $\xi_1(t),\ \xi_2(t)$.
	\end{tabular}
\end{center}

\subsubsection{The candidate potential and the role of \texorpdfstring{$P(b,c)$}{P(b,c)}}
\label{subsec:Dbc-potential}

Fix $(b,c)\in\Delta$ and write
\[
h=h_\kappa,
\qquad
m=\frac{1+b}{2}.
\]
We define $\phi=\phi_{b,c}:[0,1]\to\mathbb R$ so that the first-order condition holds along the support:
\begin{equation}
	\phi(s)
	=
	\frac12h(s^2),
	\qquad
	0\le s<b,
	\label{eq:phi-Dbc-1}
\end{equation}
\begin{equation}
	\phi(s)
	=
	\frac12h(m^2)
	+
	\int_m^s
	(1+b-y)\,
	h'\!\bigl(y(1+b-y)\bigr)\,dy,
	\qquad
	b'\le s\le c,
	\label{eq:phi-Dbc-2}
\end{equation}
\begin{equation}
	\phi(t)
	=
	\phi(c)
	+
	\int_c^t
	\frac{u\,\xi_1(u)\xi_2(u)}{\kappa}\,du,
	\qquad
	c\le t\le1,
	\label{eq:phi-Dbc-3}
\end{equation}
and
\begin{equation}
	\phi(s)
	=
	h\!\bigl(s\tau(s)\bigr)
	-
	\phi\!\bigl(\tau(s)\bigr),
	\qquad
	b\le s\le b'.
	\label{eq:phi-Dbc-4}
\end{equation}

Define the corresponding slack
\begin{equation}
	G(s,t)
	=
	\phi(s)+\phi(t)-h(st).
	\label{eq:G-Dbc}
\end{equation}

The construction is designed so that $G$ vanishes on the support.
The only non-automatic compatibility condition occurs at $s=b$.

\begin{lemma}\label{lem:potential-properties}
	Let $(b,c)\in\Delta$.
	\begin{enumerate}[label=\textup{(\alph*)}]
		\item
		Away from $b$, $\phi$ is continuously differentiable on $(0,1)$ and, for every $p\in\Gamma(t)$,
		\begin{equation}
			\phi'(t)=p\,h'(pt).
			\label{eq:first-order}
		\end{equation}
		
		\item
		\[
		G(s,t)=0
		\]
		for every $(s,t)\in\operatorname{supp}(D_{b,c})$ with $s,t\ne b$.
		
		\item
		The mismatch at $b$ is
		\begin{equation}
			\phi(b+)-\frac12h(b^2)
			=
			-\frac{P(b,c)}{24\kappa}.
			\label{eq:jump-P}
		\end{equation}
		If $b=0$, this reads
		\begin{equation}
			\phi(0)
			=
			-\frac{P(0,c)}{24\kappa}.
			\label{eq:phi-zero-P}
		\end{equation}
	\end{enumerate}
\end{lemma}

\begin{proof}
	On the diagonal and anti-diagonal pieces,
	\eqref{eq:first-order} follows directly from
	\eqref{eq:phi-Dbc-1}--\eqref{eq:phi-Dbc-2}. For $t\in[c,1]$ and
	$p\in\{\xi_1(t),\xi_2(t)\}$, the other branch is
	$q=\kappa/t-p$, so
	\[
	p\,h'(pt)
	=
	\frac{p(\kappa-pt)}{\kappa}
	=
	\frac{pqt}{\kappa}
	=
	\frac{t\,\xi_1(t)\xi_2(t)}{\kappa},
	\]
	which is exactly the derivative in \eqref{eq:phi-Dbc-3}.
	Finally, differentiating \eqref{eq:phi-Dbc-4} gives
	\[
	\phi'(s)
	=
	\tau(s)\,h'\!\bigl(s\tau(s)\bigr),
	\]
	because $s\in\Gamma(\tau(s))$.
	
	The expressions defining $\phi$ agree at $c$ by construction and at
	$\kappa/2$ because both inverse branches give
	$\tau(\kappa/2)=1$. At $b'$, agreement follows from the
	anti-diagonal identity proved below:
	\[
	\phi(b')+\phi(c)=h(b'c).
	\]
	At $b'$ and $\kappa/2$, the coupled points given by the adjacent formulas
	coincide. At $c$, the point coupled with $c$ along the anti-diagonal is $b'=\xi_2(c)$,
	and the two formulas along the curved branches give the same derivative by
	the identity above. Hence $\phi$ is continuously
	differentiable away from $b$, proving \textup{(a)}.
	
	Part \textup{(b)} is immediate on the diagonal and on the curved
	branches. On the anti-diagonal, the derivative of
	\[
	\phi(s)+\phi(1+b-s)-h\!\bigl(s(1+b-s)\bigr)
	\]
	vanishes identically, and the expression itself vanishes at
	$s=m$.
	
	For \textup{(c)}, note that $\tau(b)=c$. Hence
	\[
	\phi(b+)=h(bc)-\phi(c).
	\]
	Substituting \eqref{eq:phi-Dbc-2}, and using
	$\kappa=c(1+2b-c)$, direct integration gives
	\[
	\begin{aligned}
		24\kappa\,\phi(c)
		={}&6c^4-20c^3+18c^2-6c+1
		-32bc^3+60bc^2-24bc+4b\\
		&+42b^2c^2-30b^2c+6b^2
		-12b^3c+4b^3+b^4,
	\end{aligned}
	\]
	whereas
	\[
	24\kappa\left(h(bc)-\frac12h(b^2)\right)
	=
	24\kappa bc-12b^2c^2-12\kappa b^2+6b^4.
	\]
	Subtracting the two expressions and simplifying yields
	\[
	24\kappa\bigl(\phi(b+)-\phi(b-)\bigr)=-P(b,c),
	\]
	which proves \eqref{eq:jump-P}.
\end{proof}

Thus, for $b>0$,
\[
P(b,c)=0
\quad\Longleftrightarrow\quad
\phi_{b,c}\text{ is continuous at }b.
\]

This is the reason why the polynomial $P$ appears in the
parametrization of $\Gamma^{\mathrm{II}}$.

\subsubsection{A monotonicity lemma for the slack}
\label{subsec:slack-monotonicity}

The remaining task is to show that $G\ge0$ away from the support.
The factorization \eqref{eq:slack-factorization} turns this into a
one-dimensional monotonicity problem.

If $p=p(t)\in\Gamma(t)$, put
\[
q(t)=\frac{\kappa}{t}-p(t).
\]
Then
\begin{equation}
	\frac{\partial G}{\partial t}(s,t)
	=
	\frac{t}{\kappa}
	(s-p(t))(s-q(t)).
	\label{eq:Gt-roots}
\end{equation}
Hence, for fixed $s$, the function $t\mapsto G(s,t)$ is strictly
decreasing exactly when $s$ lies strictly between the two roots
$p(t)$ and $q(t)$.

The following lemma contains the geometric information needed below.
Local minima at the endpoints of $[0,1]$ are understood in the
one-sided sense.

\begin{lemma}\label{lem:minima-contact}
	Assume that $\phi_{b,c}$ is continuous on $[0,1]$.
	Fix $s\in[b,1]$. Every local minimum of
	\[
	t\longmapsto G(s,t)
	\]
	occurs either at $t=b$ or at a point of $\Gamma(s)$, i.e.\ on the contact set.
\end{lemma}

\begin{proof}
	We examine the two roots in \eqref{eq:Gt-roots} on the successive
	intervals determined by the support. Empty intervals are omitted;
	in particular, when $c=1$ we have $b=\kappa/2=b'$.
	
	\medskip
	\noindent\emph{The interval $0<t<b$.}
	Here $p(t)=t$. By \eqref{eq:F1},
	\[
	q(t)=\frac{\kappa}{t}-t
	>
	\frac{\kappa}{b}-b
	\ge1.
	\]
	Thus, for $s\ge b$, the function $G(s,\cdot)$ is decreasing throughout
	$(0,b)$. This interval can therefore contribute only the endpoint
	$t=b$ as a local minimum.
	
	\medskip
	\noindent\emph{The interval $b<t<\kappa/2$.}
	Here $p(t)=\tau(t)$ and $t=\xi_1(p(t))$. By
	\eqref{eq:branch-products},
	\[
	t\,p(t)\le\frac{\kappa}{2},
	\]
	and hence
	\[
	q(t)
	=
	\frac{\kappa}{t}-p(t)
	\ge
	\frac{\kappa}{2t}
	\ge1.
	\]
	Moreover $p(t)=\tau(t)$ is increasing on this interval, because it is
	the inverse of the increasing branch $\xi_1$. Consequently the set
	where $G(s,\cdot)$ decreases is either empty or an interval whose
	right endpoint is characterized by
	\[
	p(t)=s.
	\]
	At that endpoint, $t\in\Gamma(s)$. (Since $q\ge1\ge s$, with equality only at $t=\kappa/2$ and $s=1$, the decreasing set
	is exactly $\{t\in(b,\kappa/2):p(t)<s\}$; for $s=1$ it is the whole interval and its right end-point satisfies $\kappa/2\in\Gamma(1)$.)
	
	\medskip
	\noindent\emph{The interval $\kappa/2<t<1$.}
	Set
	\[
	L(t)=\min\{p(t),q(t)\},
	\qquad
	U(t)=\max\{p(t),q(t)\}.
	\]
	We first record the shape of these two functions.
	
	For $\kappa/2<t<b'$, one has $p(t)=\tau(t)$ with
	$t=\xi_2(p(t))$. Hence \eqref{eq:branch-products} gives
	$t\,p(t)\ge\kappa/2$, so that $q(t)\le p(t)$. Thus
	\[
	U(t)=p(t)=\tau(t),
	\qquad
	L(t)=q(t).
	\]
	Since $\xi_2$ is decreasing, $U$ is decreasing. To see that $L$ is
	also decreasing, write $u=p(t)$. Then
	\[
	t=\xi_2(u)
	=
	\frac{\kappa+u-u^2}{2u}
	\]
	and
	\[
	q
	=
	\frac{2\kappa u}{\kappa+u-u^2}-u.
	\]
	A direct differentiation shows that
	\[
	\frac{dq}{du}\ge0
	\quad\Longleftrightarrow\quad
	2\kappa(\kappa+u^2)
	\ge
	\bigl(\kappa+u(1-u)\bigr)^2.
	\]
	The difference between the two sides is
	\[
	\kappa^2-u^2(1-u)^2
	+
	2\kappa u(2u-1),
	\]
	which is non-negative by \eqref{eq:F4} and $u\ge c\ge1/2$.
	Since $u=\tau(t)$ decreases with $t$, $L(t)=q(t)$ is decreasing.
	
	For $b'<t<c$, $\Gamma(t)=\{p(t)\}$ with
	\[
	p(t)=1+b-t.
	\]
	The function $t\mapsto t(1+b-t)$ is concave, and its values at both
	endpoints are $b'c$. Moreover
	\[
	2b'c-\kappa=c(1-c)\ge0.
	\]
	Hence $t\,p(t)\ge\kappa/2$, so again $q(t)\le p(t)$. Therefore
	\[
	U(t)=1+b-t
	\]
	is decreasing, while
	\[
	L(t)
	=
	\frac{\kappa}{t}-(1+b-t)
	\]
	satisfies
	\[
	L''(t)=\frac{2\kappa}{t^3}>0,
	\qquad
	L(c)=b,
	\qquad
	L'(c)=1-\frac{\kappa}{c^2}\ge0
	\]
	by \eqref{eq:F3}. Thus on this interval $L$ can change monotonicity at
	most once, from decreasing to increasing.
	
	Finally, for $c<t<1$,
	\[
	L(t)=\xi_1(t),
	\qquad
	U(t)=\xi_2(t),
	\]
	so $L$ is increasing and $U$ is decreasing.
	
	The pieces match continuously at $\kappa/2$, $b'$, and $c$. We have
	therefore proved that $U$ is continuous and decreasing on
	$[\kappa/2,1]$, whereas $L$ is continuous and first decreasing and
	then increasing. In addition,
	\[
	L(c)=b,
	\qquad
	L(t)=\xi_1(t)\quad(c\le t\le1),
	\]
	and $U(t)\in\Gamma(t)$ throughout this interval.
	
	For fixed $s\ge b$, the decreasing set of $G(s,\cdot)$ in
	$(\kappa/2,1)$ is
	\[
	J_s
	=
	\{t:L(t)<s<U(t)\}.
	\]
	Since $L$ is quasi-convex and $U$ is decreasing, $J_s$ is an interval
	(possibly empty). Its right endpoint can occur only when
	\[
	U(t)=s
	\qquad\text{or}\qquad
	L(t)=s
	\]
	on the increasing branch of $L$.
	
	If $U(t)=s$, then $t\in\Gamma(s)$, because $U(t)\in\Gamma(t)$ and $D_{b,c}$ is symmetric. If $L(t)=s$, note that the increasing part of $L$
	before $c$ takes values at most $L(c)=b\le s$. Hence such a right
	endpoint is either $t=c$ with $s=b$, or lies in $[c,1]$, where
	$L=\xi_1$. In both cases $(s,t)$ belongs to the support.
	
	It remains to check the endpoints and junctions. If $b>0$, the
	function decreases immediately to the right of $0$, so $0$ is not
	a local minimum; if $b=0$, this endpoint is already allowed.
	At $t=1$ the left derivative is
	\[
	\partial_tG(s,1^-)=\frac1\kappa(s-\kappa/2)^2.
	\]
	Unless $s=\kappa/2$, it is strictly positive near $1$, excluding
	a one-sided minimum there. When $s=\kappa/2$, $1\in\Gamma(s)$.
	
	At an interior junction $t=\kappa/2>b$, both roots equal $1$.
	If $s<1$, the derivative is positive on both sides sufficiently
	close to the junction; if $s=1$, the junction itself belongs to $\Gamma(1)$.
	At $b'$ and $c$, continuity of the roots shows that a sign change
	from negative to positive is one of the right endpoints already
	classified. Coincident junctions are covered by the same one-sided
	argument. The possible junction $t=b$ is allowed separately.
	
	On each remaining open piece the derivative has only isolated
	zeros: the branches of the support are strictly monotone, and the
	complementary root is monotone or strictly convex on the pieces
	considered above. Thus there are no constant intervals that could
	produce additional local minima. Every interior local minimum
	must end an interval of decrease, and we have shown that each such
	endpoint is either $b$ or a point of $\Gamma(s)$.
\end{proof}

This lemma replaces the full two-dimensional sign analysis by a single
principle: for $s\ge b$, the slack can attain a new minimum only on the
support or at the junction $t=b$.

\subsubsection{Optimality of the candidate potential}
\label{subsec:potential-IIIII}

\begin{theorem}\label{thm:Dbc-potential}
	Let $(b,c)\in\Delta$. Assume that either
	\[
	b>0
	\quad\text{and}\quad
	P(b,c)=0,
	\]
	or
	\[
	b=0
	\quad\text{and}\quad
	P(0,c)\le0.
	\]
	Then $\phi_{b,c}$ is a Kantorovich potential for $D_{b,c}$. Hence
	\begin{equation}
		r^*\!\left(
		\frac1{12}
		+\frac14\,\E_{D_{b,c}}[S^2T^2]
		\right)
		=
		3\,\E_{D_{b,c}}[ST].
		\label{eq:Dbc-boundary}
	\end{equation}
\end{theorem}

\begin{proof}
	By Lemma~\ref{lem:potential-properties}, when $b>0$ the condition
	$P(b,c)=0$ makes $\phi$ continuous at $b$ and gives
	\[
	G(b,b)=0.
	\]
	When $b=0$, the condition $P(0,c)\le0$ gives
	\[
	G(0,0)=2\phi(0)\ge0.
	\]
	
	Under the assumptions of the theorem, $\phi$ is continuous on
	$[0,1]$ (at the endpoint $0$ when $b=0$, no matching condition is
	needed). Hence, for fixed $s$, the function $G(s,\cdot)$ is continuous
	on $[0,1]$ and piecewise continuously differentiable.
	
	We first prove
	\[
	G(s,t)\ge0
	\qquad
	\text{for }s\in[b,1].
	\]
	Take first $s=b$. By Lemma~\ref{lem:minima-contact}, every local
	minimum of $G(b,\cdot)$ occurs either at $t=b$ or at a point of $\Gamma(b)$.
	If $b>0$, then $G(b,b)=0$ because $P(b,c)=0$. If $b=0$, we already
	know that $G(0,0)\ge0$. The other relevant point of $\Gamma(b)$ is $t=c$ (the point $(b,c)=(\xi_1(c),c)$ lies on the curved branch), and
	\eqref{eq:phi-Dbc-4} gives directly
	\[
	G(b,c)=0.
	\]
	Thus every local minimum of $G(b,\cdot)$ is non-negative. Since the
	function is continuous on the compact interval $[0,1]$,
	\[
	G(b,t)\ge0
	\qquad
	0\le t\le1.
	\]
	By symmetry,
	\[
	G(s,b)\ge0
	\qquad
	0\le s\le1.
	\]
	
	Now fix any $s\in[b,1]$. Again by
	Lemma~\ref{lem:minima-contact}, every local minimum of
	$G(s,\cdot)$ occurs either at a point of $\Gamma(s)$, where $G=0$, or at
	$t=b$, where $G(s,b)\ge0$. Therefore
	\[
	G(s,t)\ge0
	\qquad
	b\le s\le1,\quad 0\le t\le1.
	\]
	
	It remains to consider $0\le s<b$, which occurs only when $b>0$.
	For $t<b$,
	\[
	G(s,t)
	=
	\frac{(s-t)^2}{2}
	\left(
	1-\frac{(s+t)^2}{2\kappa}
	\right)
	\ge0,
	\]
	because $4b^2\le2\kappa$ by \eqref{eq:F1} (indeed $\kappa\ge b(1+b)\ge2b^2$).
	
	Now fix $t\ge b$ and let
	\[
	F(s)=G(s,t),
	\qquad
	0\le s\le b.
	\]
	Then
	\[
	F''(s)
	=
	1-\frac{3s^2}{\kappa}+\frac{t^2}{\kappa}
	\ge
	1-\frac{2b^2}{\kappa}
	\ge0,
	\]
	while
	\[
	F'(b^-)
	=
	(b-t)
	\left(
	1-\frac{b(b+t)}{\kappa}
	\right)
	\le0
	\]
	by \eqref{eq:F1}. Since $F$ is convex, $F'$ is non-decreasing, so $F'(s)\le F'(b^-)\le0$ on $[0,b]$. Hence $F$ is non-increasing on $[0,b]$, and
	\[
	G(s,t)=F(s)\ge F(b)=G(b,t)\ge0.
	\]
	
	Thus $G\ge0$ on $[0,1]^2$ and $G=0$ on
	$\operatorname{supp}(D_{b,c})$. For $b>0$ the diagonal junction
	$(b,b)$ satisfies equality by the matching condition; for $b=0$
	the point $(0,0)$ is absent from the support. Therefore $\phi_{b,c}$ is a Kantorovich
	potential, and \eqref{eq:Dbc-boundary} follows from
	Lemma~\ref{lem:kantorovich-sufficient}.
\end{proof}

\subsubsection{Selection of the boundary parameters}
\label{subsec:parameter-selection}

It remains to determine when the assumptions of
Theorem~\ref{thm:Dbc-potential} hold.

\begin{lemma}\label{lem:P-selection}
	For
	\[
	0\le b\le\frac15,
	\qquad
	\frac12+b\le c\le1,
	\]
	the following statements hold.
	\begin{enumerate}[label=\textup{(\roman*)}]
		\item
		In the interior of $\Delta$,
		\begin{equation}
			\frac{\partial P}{\partial c}(b,c)
			=
			6(2c-1-2b)
			\left(
			2(1-c)^2-(1-b)^2
			\right)
			<0.
			\label{eq:Pc-negative}
		\end{equation}
		
		\item
		\begin{equation}
			P\!\left(b,\frac12+b\right)
			=
			\frac18
			\left(
			3-12b^2+16b^3-8b^4
			\right)
			>0,
			\label{eq:P-left-endpoint}
		\end{equation}
		whereas
		\begin{equation}
			P(b,1)
			=
			(1-b)^3(5b-1).
			\label{eq:P-right-endpoint}
		\end{equation}
		
		\item
		For every $b\in[0,1/5]$ there is a unique
		\[
		c(b)\in\left(\frac12+b,1\right]
		\]
		such that $P(b,c(b))=0$. Moreover,
		\[
		c(0)=c^*,
		\qquad
		c\!\left(\frac15\right)=1,
		\]
		and
		\[
		P(0,c)\le0
		\quad\Longleftrightarrow\quad
		c^*\le c\le1.
		\]
	\end{enumerate}
\end{lemma}

\begin{proof}
	The identities in \textup{(i)} and \textup{(ii)} follow by direct
	factorization. In the interior of $\Delta$,
	\[
	2c-1-2b>0,
	\]
	while
	\[
	2(1-c)^2
	\le
	\frac12(1-2b)^2
	<
	(1-b)^2.
	\]
	Hence $\partial_cP<0$.
	
	For $b\le\frac15$ we have $16b^3-8b^4=8b^3(2-b)\ge0$ and $12b^2\le\frac{12}{25}<3$, so $P(b,\frac12+b)>0$, while
	$P(b,1)\le0$ with equality only for $b=\frac15$. The endpoint values in \textup{(ii)}, together with the strict
	monotonicity in $c$, give existence and uniqueness of the zero.
	The formula for $\partial_cP$ is strictly negative along the
	entire zero curve, including $b=0$ and $b=1/5$.
	The implicit function theorem therefore makes $b\mapsto c(b)$
	continuous and real analytic locally, with analytic extensions
	across its two endpoints. The remaining statements follow by evaluating the endpoint
	cases $b=0$ and $b=1/5$.
\end{proof}

\begin{corollary}\label{cor:familiesIIIII-boundary}
	For every $c\in[c^*,1]$,
	\[
	r^*\!\left(\ARV(C^{\mathrm{III}}_c)\right)
	=
	\rho(C^{\mathrm{III}}_c),
	\]
	and for every $b\in[0,1/5]$,
	\[
	r^*\!\left(\ARV(C^{\mathrm{II}}_{b,c(b)})\right)
	=
	\rho(C^{\mathrm{II}}_{b,c(b)}),
	\]
	where at $b=0$ we use
	\[
	C^{\mathrm{II}}_{0,c^*}=C^{\mathrm{III}}_{c^*}.
	\]
	Consequently, $\Gamma^{\mathrm{III}}$ and
	$\Gamma^{\mathrm{II}}$ are upper-boundary arcs.
\end{corollary}

\begin{proof}
	If $b=0$ and $c\in[c^*,1)$, then $P(0,c)\le0$, so
	Theorem~\ref{thm:Dbc-potential} applies to $D_{0,c}$. The endpoint $c=1$ is
	$C^{\mathrm{III}}_1=C_2$. By Proposition~\ref{prop:region-geometry},
	$r^*(11/120)=1/2$, and $C_2$ attains this value. Hence this point
	belongs to the upper boundary, as required.
	
	If $0<b\le1/5$, then $P(b,c(b))=0$, and the same theorem applies to
	$D_{b,c(b)}$. The case $b=0$ is the common endpoint of the two arcs.
\end{proof}

\subsubsection{Completion of the proof of the main theorem}
\label{subsec:main-proof}

\begin{proof}[Proof of Theorem~\ref{thm:main}]
	Theorem~\ref{thm:familyI-optimal} shows that
	$\Gamma^{\mathrm I}$ lies on the upper boundary, while
	Corollary~\ref{cor:familiesIIIII-boundary} gives the same conclusion
	for $\Gamma^{\mathrm{II}}$ and $\Gamma^{\mathrm{III}}$.
	
	The three arcs join continuously because
	\[
	C^{\mathrm{II}}_{0,c^*}
	=
	C^{\mathrm{III}}_{c^*},
	\qquad
	C^{\mathrm{II}}_{1/5,1}
	=
	C^{\mathrm I}_{2/5}.
	\]
	Moreover,
	\[
	C^{\mathrm{III}}_1=C_2,
	\qquad
	C^{\mathrm I}_0=M.
	\]
	Hence their union is a continuous subset of the upper boundary joining
	\[
	\left(\frac{11}{120},\frac12\right)
	\quad\text{to}\quad
	\left(\frac{2}{15},1\right).
	\]
	Its first coordinate therefore takes every value in
	$[11/120,2/15]$. Since the upper boundary is the graph of the function
	$r^*$ defined in \eqref{eq:rstar-folded}, we conclude that
	\[
	\operatorname{graph}(r^*)
	=
	\Gamma^{\mathrm{III}}
	\cup
	\Gamma^{\mathrm{II}}
	\cup
	\Gamma^{\mathrm I}.
	\]
	
	Proposition~\ref{prop:region-geometry} supplies the description of
	the full region, its symmetry, and the continuity, concavity and
	strict monotonicity of $r^*$. The first and third arcs have analytic
	parametrizations by their displayed moment formulas, while the
	analyticity of the middle parametrization follows from
	Lemma~\ref{lem:P-selection}. This completes the proof.
\end{proof}

\subsubsection{Regularity of the boundary and behavior near the minimum}
\label{subsec:vertical-tangent}

We now prove \cref{cor:C1,cor:vertical-tangent}. The optimality proofs
provide, at every point of the three arcs, a supporting line with explicit
slope, which yields the regularity of the boundary.

\begin{proof}[Proof of \cref{cor:C1}]
	\emph{Step 1: supporting lines.} Let $D$ be the folded copula of a member
	of one of the three families, with associated $\kappa>0$ (that is,
	$\kappa=\kappa_b$ for Family~I, and $\kappa=c(1+2b-c)$ for
	Families~II and~III, with $c<1$ when $b=0$). By
	Theorem~\ref{thm:familyI-optimal} and
	Theorem~\ref{thm:Dbc-potential}, $D$ maximizes
	$\E_{D'}[c_\kappa(S,T)]=\E_{D'}[ST]-\E_{D'}[S^2T^2]/(2\kappa)$ over
	$D'\in\mathcal C$. By \eqref{eq:fold-arv}, \eqref{eq:fold-rho} and
	$\rho(\unfold D)=3\E_D[ST]$, for every copula $C$,
	\begin{align*}
		\rho(C)-\frac6\kappa\ARV(C)
		&\le3\,\E_{\fold C}[c_\kappa(S,T)]-\frac1{2\kappa}
		\le3\,\E_{D}[c_\kappa(S,T)]-\frac1{2\kappa}\\
		&=\rho(\unfold D)-\frac6\kappa\ARV(\unfold D).
	\end{align*}
	Hence the line of slope $6/\kappa$ through the corresponding boundary
	point supports $\A$.
	
	\emph{Step 2: slopes.} For $a=0$ ($D=M$, $\kappa=2$), Step~1 gives
	$r^*(\alpha)\le1+3(\alpha-2/15)$; hence the left derivative of $r^*$
	at $2/15$ is at least $3$ and, by concavity, every supergradient of
	$r^*$ at an interior point is at least $3$.
	
	\emph{Step 3: no corners.} Suppose that, for some
	$\alpha_0\in(11/120,2/15)$, the one-sided derivatives satisfy
	$(r^*)'_+(\alpha_0)<(r^*)'_-(\alpha_0)$. Choose a nonempty open interval
	$J$ of slopes strictly between them. Shrinking $J$ if necessary, we may
	assume that $K=\{6/m:m\in J\}$ is contained in one of the intervals
	$(0,\kappa^*)$, $(\kappa^*,\frac25)$, $(\frac25,2)$, where
	$\kappa^*=c^*(1-c^*)$. For $m\in J$, concavity shows that the line of
	slope $m$ through $(\alpha_0,r^*(\alpha_0))$ meets $\A$ only at this
	point. Along the parameter intervals $c\in(c^*,1)$ (Family~III),
	$b\in(0,\frac15)$ (Family~II, with $c=c(b)$) and $a\in(0,\frac25)$
	(Family~I), $\kappa$ depends continuously on the parameter and tends to
	the endpoint values $0,\kappa^*$, respectively $\kappa^*,\frac25$, respectively
	$\frac25,2$. By the intermediate value theorem, the parameters of the
	corresponding family with $\kappa\in K$ form a nonempty open set,
	which contains an open interval $I$. By Step~1 every parameter in $I$
	generates the point $(\alpha_0,r^*(\alpha_0))$, so $\ARV$ is constant
	on $I$. This is impossible: along each family, $\ARV$ is a real-analytic
	function of the parameter (by \eqref{eq:familyI-arv-rho},
	\eqref{eq:familyII-arv}, \eqref{eq:familyIII-arv} and
	Lemma~\ref{lem:P-selection}) that takes different values at the two
	endpoints of the parameter interval, hence it is not constant on any
	open subinterval.
	
	\emph{Step 4: conclusion.} By Step~3, $r^*$ is differentiable on
	$(11/120,2/15)$; a concave function that is differentiable on an open
	interval is continuously differentiable there. At a point generated by a
	family member with parameter $\kappa$, the supporting line of Step~1
	must then be the tangent line, which gives \eqref{eq:rstar-derivative};
	for Family~I, $\kappa_{1-2a}=\frac23(1-a)(3-5a)$. As $a\downarrow0$ the
	derivative tends to $3$, which is therefore the left derivative of the
	concave function $r^*$ at $2/15$. The values at the junctions follow
	from $\kappa_{1/5}=\frac25$ and $c(0)=c^*$. Finally, as $c\uparrow1$,
	$\ARV(C^{\mathrm{III}}_c)\downarrow11/120$ and $6/(c(1-c))\to+\infty$.
\end{proof}

\begin{proof}[Proof of \cref{cor:vertical-tangent}]
	Near the left endpoint the upper boundary is generated by
	$C^{\mathrm{III}}_c$: by the expansion below, $c\mapsto\ARV(C^{\mathrm{III}}_c)$ is continuous and strictly decreasing near $c=1$,
	with value $\frac{11}{120}$ at $c=1$, so every $\alpha$ close to $\frac{11}{120}$ equals $\ARV(C^{\mathrm{III}}_c)$ for a unique $c$ close to $1$, and then
	$r^*(\alpha)=\rho(C^{\mathrm{III}}_c)$. Set $c=1-\delta$, with $\delta\downarrow0$.
	Expanding \eqref{eq:familyIII-rho} and \eqref{eq:familyIII-arv} gives
	\begin{align*}
		\rho(C^{\mathrm{III}}_{1-\delta})
		&=\frac12+\frac{\delta^3}{2}-\frac{\delta^4}{2}+O(\delta^5),\\
		\ARV(C^{\mathrm{III}}_{1-\delta})
		&=\frac{11}{120}+\frac{\delta^4}{16}-\frac{7\delta^5}{60}+O(\delta^6).
	\end{align*}
	Hence
	\[
	\delta
	=
	2\left(\alpha-\frac{11}{120}\right)^{1/4}
	+
	o\!\left(
	\left(\alpha-\frac{11}{120}\right)^{1/4}
	\right),
	\]
	and substitution into the expansion of $\rho$ yields the stated
	formula.
\end{proof}


\begin{thebibliography}{99}

\bibitem{AndroschuckEtAl2017}
  T.~Androschuck et al.,
  Simulation-based capital models: testing, justifying and communicating
  choices. A report from the Life Aggregation and Simulation Techniques
  Working Party, \emph{British Actuarial Journal} 22(2) (2017) 257--335.

\bibitem{AnsariRockel2026}
  J.~Ansari, M.~Rockel,
  The exact region and an inequality between Chatterjee's and Spearman's rank correlations,
  \emph{J.\ Multivariate Anal.} 214 (2026), 105630.

\bibitem{AnzilliDurante2026a}
  L.~Anzilli, F.~Durante,
  The average rectangular volume induced by supermodular aggregation functions,
  \emph{J.\ Math.\ Anal.\ Appl.} 555 (2026), 130137.

\bibitem{AnzilliDurante2026b}
  L.~Anzilli, F.~Durante,
  On the $(\ARV,\rho)$ region for copulas,
  in: \emph{Soft Methods in Probability and Statistics for Complex Data}, Proceedings of SMPS 2026, Springer, to appear.

\bibitem{CuadrasAuge1981}
  C.\,M.~Cuadras, J.~Aug\'e,
  A continuous general multivariate distribution and its properties,
  \emph{Comm.\ Statist.\ Theory Methods} 10(4) (1981) 339--353.

\bibitem{Daniels1950}
  H.\,E.~Daniels,
  Rank correlation and population models,
  \emph{J.\ R.\ Stat.\ Soc.\ Ser.\ B} 12 (1950) 171--181.

\bibitem{DuranteKolesarovaMesiarSempi2008}
  F.~Durante, A.~Koles\'arov\'a, R.~Mesiar, C.~Sempi,
  Semilinear copulas,
  \emph{Fuzzy Sets Syst.} 159(1) (2008) 63--76.

\bibitem{DuranteSarkociSempi2009}
  F.~Durante, P.~Sarkoci, C.~Sempi,
  Shuffles of copulas,
  \emph{J.\ Math.\ Anal.\ Appl.} 352 (2009) 914--921.

\bibitem{DuranteSempi2016}
  F.~Durante, C.~Sempi,
  \emph{Principles of Copula Theory},
  CRC/Chapman \& Hall, Boca Raton, 2016.

\bibitem{FuchsLimbachSchurrer2026}
  S.~Fuchs, C.~Limbach, F.~Sch\"urrer,
  On exact regions between measures of concordance and Chatterjee's rank correlation for lower semilinear copulas,
  \emph{Internat.\ J.\ Approx.\ Reason.} 189 (2026), 109588.

\bibitem{GangboMcCann1996}
  W.~Gangbo, R.\,J.~McCann,
  The geometry of optimal transportation,
  \emph{Acta Math.} 177 (1996) 113--161.

\bibitem{GrabischMarichalMesiarPap2009}
  M.~Grabisch, J.-L.~Marichal, R.~Mesiar, E.~Pap,
  \emph{Aggregation Functions}, Cambridge University Press, Cambridge, 2009.

\bibitem{KlementMesiarPap2010}
  E.\,P.~Klement, R.~Mesiar, E.~Pap,
  A universal integral as common frame for Choquet and Sugeno integral,
  \emph{IEEE Trans.\ Fuzzy Syst.} 18 (2010) 178--187.

\bibitem{KokolBukovsekEtAl2021}
  D.~Kokol Bukov\v{s}ek, T.~Ko\v{s}ir, B.~Moj\v{s}kerc, M.~Omladi\v{c},
  Spearman's footrule and Gini's gamma: local bounds for bivariate copulas
  and the exact region with respect to Blomqvist's beta,
  \emph{J.\ Comput.\ Appl.\ Math.} 390 (2021), 113385.

\bibitem{KokolBukovsekEtAl2026}
  D.~Kokol Bukov\v{s}ek, P.~Lazi\'c, B.~Moj\v{s}kerc, N.~Stopar,
  The exact region determined by Spearman's footrule, Gini's gamma and Kendall's tau,
  \emph{J.\ Comput.\ Appl.\ Math.} 486 (2026), 117639.

\bibitem{KokolBukovsekMojskerc2022}
  D.~Kokol Bukov\v{s}ek, B.~Moj\v{s}kerc,
  On the exact region determined by Spearman's footrule and Gini's gamma,
  \emph{J.\ Comput.\ Appl.\ Math.} 410 (2022), 114212.

\bibitem{KokolBukovsekStopar2023}
  D.~Kokol Bukov\v{s}ek, N.~Stopar,
  On the exact regions determined by Kendall's tau and other concordance measures,
  \emph{Mediterr.\ J.\ Math.} 20(3) (2023), 147.

\bibitem{KokolBukovsekStopar2024}
  D.~Kokol Bukov\v{s}ek, N.~Stopar,
  On the exact region determined by Spearman's rho and Spearman's footrule,
  \emph{J.\ Comput.\ Appl.\ Math.} 437 (2024), 115463.

\bibitem{Lancaster1958}
  H.\,O.~Lancaster,
  The structure of bivariate distributions,
  \emph{Ann.\ Math.\ Statist.} 29(3) (1958) 719--736.

\bibitem{McNeilNeslehovaSmith2025}
  A.\,J.~McNeil, J.\,G.~Ne\v{s}lehov\'a, A.\,D.~Smith,
  Measures and models of non-monotonic dependence,
  arXiv:2512.10828, 2025.

\bibitem{MikusinskiSherwoodTaylor1992}
  P.~Mikusi\'nski, H.~Sherwood, M.\,D.~Taylor,
  Shuffles of Min,
  \emph{Stochastica} 13 (1992) 61--74.

\bibitem{Nelsen2006}
  R.\,B.~Nelsen,
  \emph{An Introduction to Copulas}, 2nd ed.,
  Springer, New York, 2006.

\bibitem{RachevRuschendorf1998}
  S.\,T.~Rachev, L.~R\"uschendorf,
  \emph{Mass Transportation Problems. Vol.~I: Theory},
  Springer, New York, 1998.

\bibitem{Rockel2026}
  M.~Rockel,
  The exact region between Chatterjee's and Blest's rank correlations,
  \emph{Internat.\ J.\ Approx.\ Reason.} 197 (2026), 109744.

\bibitem{Santambrogio2015}
  F.~Santambrogio,
  \emph{Optimal Transport for Applied Mathematicians},
  Progress in Nonlinear Differential Equations and Their Applications 87,
  Birkh\"auser, Cham, 2015.

\bibitem{SchreyerPaulinTrutschnig2017}
  M.~Schreyer, R.~Paulin, W.~Trutschnig,
  On the exact region determined by Kendall's $\tau$ and Spearman's $\rho$,
  \emph{J.\ R.\ Stat.\ Soc.\ Ser.\ B} 79 (2017) 613--633.

\bibitem{ShawSmithSpivak2011}
  R.\,A.~Shaw, A.\,D.~Smith, G.\,S.~Spivak,
  Measurement and modelling of dependencies in economic capital,
  \emph{British Actuarial Journal} 16(3) (2011) 601--699.

\bibitem{Tschimpke2025}
  M.~Tschimpke,
  On the exact region determined by Spearman's $\rho$ and Blest's measure
  of rank correlation $\nu$ for bivariate extreme-value copulas,
  \emph{J.\ Multivariate Anal.} 205 (2025), 105377.

\bibitem{TschimpkeSchreyerTrutschnig2025}
  M.~Tschimpke, M.~Schreyer, W.~Trutschnig,
  Revisiting the region determined by Spearman's $\rho$ and Spearman's footrule $\phi$,
  \emph{J.\ Comput.\ Appl.\ Math.} 457 (2025), 116259.

\bibitem{Villani2009}
  C.~Villani,
  \emph{Optimal Transport. Old and New},
  Grundlehren der mathematischen Wissenschaften 338, Springer, Berlin, 2009.

\end{thebibliography}
\end{document}